\documentclass[a4paper,11pt]{amsart}

\usepackage{amsmath}
\usepackage{amssymb}
\usepackage{amstext}
\usepackage{amsbsy}
\usepackage{amsopn}
\usepackage{amsthm}
\usepackage{amscd}
\usepackage{amsxtra}
\usepackage{amsfonts}

\usepackage{ifthen}
\usepackage{xspace}
\usepackage{layout}
\usepackage{hyperref}

\input xy
\xyoption{matrix}
\xyoption{arrow}
\xyoption{curve}
\newdir{ (}{{}*!/-5pt/@^{(}}

\newcommand{\xycenter}[1]{\begin{center}
                          \mbox{\xymatrix{#1}}
                          \end{center}
                         }

\newboolean{xlabels}
\newcommand{\xlabel}[1]{
                        \label{#1}
                        \ifthenelse{\boolean{xlabels}}
                                   {\marginpar{{\tiny #1}}}
                                   {}
                       }

\newcommand{\ZZ}{\mathbb{Z}}
\newcommand{\NZ}{\mathbb{N}}

\newcommand{\pz}{\mathbb{P}}

\newcommand{\CC}{\mathbb{C}}

\newcommand{\PP}{\mathbb{P}}
\newcommand{\GG}{\mathbb{G}}
\newcommand{\FF}{\mathbb{F}}

\newcommand{\sF}{{\mathcal F}}

\newcommand{\sI}{{\mathcal I}}

\newcommand{\sN}{{\mathcal N}}
\newcommand{\sO}{{\mathcal O}}

\newcommand{\sQ}{{\mathcal Q}}

\newcommand{\sU}{{\mathcal U}}

\newcommand{\suchthat}{\, | \,}

\newboolean{probleme}
\newcommand{\problem}[1]
           {\ifthenelse{\boolean{probleme}}
                       {{\bf(PROBLEM: #1)\bf}}
                       {}
           }
\newboolean{zukuenftiges}
\newcommand{\zukunft}[1]
           {\ifthenelse{\boolean{zukuenftiges}}
                       {{\bf(AUSBAUM\"OGLICHKEIT: #1)\bf}}
                       {}
           }
\newboolean{extras}
\newcommand{\extra}[1]
           {\ifthenelse{\boolean{extras}}
                       {{\bf EXTRA #1 EXTRA\bf}}
                       {}
           }

\newboolean{ignore}
\newcommand{\ignore}[1]
           {\ifthenelse{\boolean{ignore}}
                       {{\bf IGNORE #1 IGNORE\bf}}
                       {}
           }

\newcommand{\longto}[1]{\stackrel{#1}{\longrightarrow}}

\DeclareMathOperator{\codim}{codim}
\DeclareMathOperator{\coker}{coker}

\DeclareMathOperator{\Aut}{Aut}

\DeclareMathOperator{\Hom}{Hom}

\DeclareMathOperator{\rank}{rank}

\DeclareMathOperator{\supp}{supp}

\DeclareMathOperator{\gin}{gin}

\theoremstyle{plain}
\begingroup

\newtheorem{thm}{Theorem}
\newtheorem{cor}[thm]{Corollary}
\newtheorem{lem}[thm]{Lemma}

\newtheorem{prop}[thm]{Proposition}

\numberwithin{thm}{subsection} 
\endgroup

\newtheorem*{thm*}{Theorem}
\newtheorem*{conj*}{Conjecture}
\newtheorem*{verm*}{Vermutung}

\theoremstyle{definition}
\newtheorem{defn}[thm]{Definition}
\newtheorem{rem}[thm]{Remark}
\newtheorem{example}[thm]{Example}

\newtheorem{notation}[thm]{Notation}

\newtheorem{construction}[thm]{Construction}

\numberwithin{equation}{section}

\newcommand{\nosubsections}{\renewcommand{\thethm}{\thesection.\arabic{thm}}

\setcounter{thm}{0}
                           }

\newcommand{\cref}[3]{(\ref{#1}, #2 \ref{#3})}

\setboolean{probleme}{false}
\setboolean{xlabels}{false}

\begin{document}
\title{Classification of rational surfaces of degree 11 and sectional 
genus 11 in $\PP^4$}
\author[v. Bothmer]{Hans-Christian Graf v. Bothmer}
\address{Institiut f\"ur algebraische Geometrie, Leibnitz Universit\"at Hannover, 
Welfengarten 
1, D-30167 Hannover }
\urladdr{\href{http://www-ifm.math.uni-hannover.de/~bothmer}{http://www-ifm.math.uni-hannover.de/\textasciitilde 
bothmer}}
\author[Ranestad]{Kristian Ranestad}
\thanks{Supported by the Schwerpunktprogramm ``Global Methods in 
Complex
        Geometry'' of the Deutsche Forschungs Gemeinschaft and by the 
strategic university program Suprema of  NFR (project
        154077/420)}
\urladdr{\href{http://folk.uio.no/ranestad/}{http://folk.uio.no/ranestad/}}

\begin{abstract}
We use the BGG-correspondence to show that there are at most three
possible Hilbert functions for smooth rational surfaces of degree 11 
and
sectional genus 11. Surfaces with one of these  Hilbert functions 
have been classified
by Popescu. The classification for a second one is done in this 
paper. For the
third Hilbert function the classification is still open.
\end{abstract}

\maketitle

\section{Introduction}
\nosubsections

In the classification of smooth embedded varieties a natural start is to determine which Hilbert polynomials occur. Next 
one can classify the
possible Hilbert functions.  A bold aim is to
determine the irreducible components of the 
Hilbert scheme
representing smooth varieties and give a description of a general 
member in each component.  

For space 
curves the first question was answered by 
Gruson and Peskine \cite {GPgenre2}, while the second and third are 
only
partially answered.
For smooth surfaces in $\PP^4$ there are even fewer results.  The 
only general 
ones are the bounds of Ellingsrud and Peskine \cite 
{EllingsrudPeskine} that give asymptotic restrictions on the Hilbert 
polynomials that
occur.   Most work has concentrated on 
small invariants, and even here the results are only partial:  The 
first question is only completely answered up to degree 10, on the
other hand in this case also the second and third question is
completely answered, allthough not explicitly.

The techniques involved in the classification of surfaces with small
invariants have developed considerably over time.  The only common 
feature is that the combination of different approaches required 
often give the arguments a certain ad hoc flavour.   While proving that 
certain components are nonempty, i.e. to give examples, can often be 
done
transparently, it is the converse result, that a component is empty,
that seems to require a combination of techniques.

In this paper we use the Tate resolution of the ideal sheaf $\sI:=\sI_{S}$ of 
a smooth surface $S$
in $\PP^4$ with Hilbert polynomial $P_{S}(n)=\frac {11}2 n^2-\frac 92 
n+1$ to determine
the Hilbert function  of $S$.
The main new idea is to study complexes on the Grassmannian of linear 
subspaces
of $\PP^4$, as introduced by \cite{ESchow}, that are analogs of the Beilinson Monad.
The degeneracy loci of the maps of these complexes define special 
linear subspaces in $\PP^4$.

Special linear subspaces are those where 
the cohomology
of the restricted ideal sheaf differs from the cohomology of a 
general restriction.
Lines in the
surface and lines that intersect the surface $S$ in a scheme of large
length are special.  A special plane intersects the 
surface $S$ in a 
curve, or in $11$ points in special position.
A special hyperplane intersects the surface in a space curve  
which lies on more surfaces of small degree than the general one.  
The geometry of special linear subspaces allows us to determine which 
maps can occur in the Tate resolution.

The main results and the organization of the paper is as follows:
In Section \ref{sPrelim} we recall the basic facts about the Tate
resolution of $\sI (n)$, its relation to the Beilinson Monad and the corresponding
complexes on the Grassmannians as explained by Eisenbud and Schreyer
\cite {ESchow}. When the 
intersection between the linear subspace and the surface is improper,
then the ideal sheaf of the intersection is not the restriction 
of the ideal sheaf. The difference is made precise by the excess conomormal sheaf.
(cf. Fulton \cite {Fu}). In Sections \ref{sLines} and \ref{sPlanes} we use this excess sheaf to study the 
restriction of the
ideal of the surface to special lines and planes respectively. 
In the Section \ref{sDiagrams} we recall how the diagrams of
generic initial ideals, as introduced by Green \cite {GreenGin}, can 
be applied to 
classify plane sections of $S$.
Section \ref{sElevenEleven} is then
devoted to determining the Hilbert function of $S$.
We show that only three different functions are possible. They 
differ from the Hilbert polynomial only in degrees $n=1,\ldots,5$, 
where
their values are $(5,15,35,70,116)$, $(5,15,35,69,116)$ and 
$(5,15,35,69,115)$
respectively.
Popescu showed that the first function occurs with three different
irreducible families. In \cite{newfamily} v.\,Bothmer, Erdenberger 
and Ludwig 
give an example with the second function which was found by a random 
search over $\FF_{2}$. 
In Section \ref{sa1b0} and \ref{sDimensions} we show that all smooth surfaces 
with this second Hilbert function belong to the same
irreducible and unirational family. This is Theorem \ref{tcons}.
In Section \ref{sConstruction} we a geometric construction of the surfaces in this 
family  
(Theorem \ref{tlinkageconstruction}).

The third Hilbert function also occurs for an irreducible
family of surfaces, but we are not able to determine whether any of 
the
surfaces belonging to that family are smooth. This reflects the 
nature of our 
methods.  The exterior algebra methods we employ do not
distinguish between smooth and singular irreducible surfaces.  It is 
in 
combination with geometric arguments that we are sometimes able to
make that distinction.   On the other hand the constructed examples 
needed to eventually prove that a component is nonempty are often so
rigid that Bertini type theorems do not easily apply. 
Therefore, explicit examples of smooth surfaces are constructed using 
the computer algebra program Macaulay2 \cite {M2}. Scripts are 
provided and 
documented on our website \cite{linesplanesweb}.  These 
examples are constructed algebraically over the $\ZZ$, and computed 
over a finite field, so by the openness condition of smoothness they 
are smooth over the rational numbers, and hence also over $\CC$. (cf. 
\cite{DES}, Appendix A).

\section{Preliminaries} \xlabel{sPrelim}
\nosubsections

\newcommand{\fU}{U}
\newcommand{\fR}{R}
\newcommand{\sOPone}{\sO_{\PP^1}}

\begin{notation}
$\quad $

\begin{tabular}{ll}
$W$ & a vector space of dimension 5\\
$E = \bigwedge W^*$ & the exterior algebra over its dual space \\
$\PP^4 = \PP(W)$ & the Grothendiek projectivisation of $W$\\
$\GG_l$ & the Grassmannian of codim $l$ linear subspaces in $\PP^4$\\
$\FF_l$ & the Flag variety of points in codim $l$ linear subspaces of
$\PP^4$\\
\end{tabular}
\end{notation}

\newcommand{\Floystad}{Fl{\o}ystad}
\newcommand{\BGG}{Bern{\v{s}}te{\u\i}n, Gel{$'$}fand and Gel{$'$}fand}
\newcommand{\Beilinson}{Be{\u\i}linson}

In this paper we use the BGG-correspondence of  \BGG\, \cite{BGG} in 
an
explicit version described by Eisenbud, \Floystad\, and Schreyer in 
\cite{EFSsheaf}.
For every sheaf $\sF$ on $\PP(W)$ one can construct a canonical exact 
complex $T(\sF)$
over the exterior algebra $E$. This complex is called the {\sl Tate 
resolution}, see \cite[Section 4]{EFSsheaf} for the construction. The 
terms of the Tate resolution can be explicitly
described:

\begin{thm}[\BGG; Eisenbud, \Floystad\, and Schreyer] \xlabel{tTate}
If $\sF$ is a coherent sheaf on $\PP(W)$, then the $e$-th term of the 
Tate resolution
is
\[
T(\sF)^e = \bigoplus_{j} \Hom_K(E,H^j(\sF(e-j))).
\]
\end{thm}

\begin{proof}
\cite[Theorem 4.1]{EFSsheaf}
\end{proof}

Now consider the incidence correspondence
\xycenter{ \FF_l \ar[d]^{\pi_1} \ar[r]^{\pi_2}& \GG_l \\
\PP^n}
and the tautological sequence
\[
 0 \to \sU_l \to W\otimes \sO_{\GG_l} \to \sQ_l \to 0.
\]
In \cite{ESchow} Eisenbud and Schreyer define an additive functor 
$\fU_l$
from graded free modules over $E$ to locally free sheaves on $\GG_l$ 
by
taking $\fU_l(E(p)) = \bigwedge^p \sU_l$ and sending a map
$\eta \colon E(q) \to E(q-p)$ to the map
$\fU_l(\eta) \colon \bigwedge^q \sU_l \to \bigwedge^{q-p}\sU_l$ 
defined by
the element of $\bigwedge^p W^*$ corresponding to $\eta$ 
\cite[Proposition 1.1]{ESchow}. 
We write $\fU_l(\sF)$
for $\fU_l(T(\sF))$.

\begin{thm}[\Beilinson; Eisenbud and Schreyer]
If $\sF$ is a sheaf on $\PP^n$ then
\[
\fU_l(\sF) \sim \fR{\pi_2}_*(\pi_1^*\sF)
\]
in the derived category.
\end{thm}

\begin{proof} \cite[Theorem 1.2]{ESchow}. \end{proof}

\begin{rem}
Notice that $\fU_n(\sF)$ is the \Beilinson-Monad \cite{Beilinson}. In 
this case $\sF$
appears as the homology
in step $0$. For $l \ge \dim \supp \sF$ one can recover $\sF$ from
$\fU_l(\sF)$. \cite[Proposition 1.3]{ESchow}
In this paper we also use the partial information contained in
$\fU_l(\sF)$ for $l \le \dim \supp \sF$.
\end{rem}

We now consider the twisted ideal sheaf $\sI_S(n)$ of a smooth
surface in $\PP^4$ and want to determine its Tate resolution 
$T(\sI_S(n))$. Its terms
\[
T(\sI_S(n))_e = \sum_{j=0}^4 H^j(\sI_S(n+e-j))\otimes E(j-e)
\]
are given by Theorem \ref{tTate}. To determine its maps
we apply
the functor $\fU_l$ and use the geometry of $S$ to analyse the 
complexes
$\fU_l(\sF)$. For fixed $l$ we use the notation
\[
F_e := \fU_l(\sF)_e = \sum_{j=0}^4 H^j(\sI_S(n+e-j))\otimes 
\bigwedge^{j-e}\sU_l.
\]
Notice that only the terms $F_{-l} ... F_{4}$ are nonzero, so the
complex $\fU_l(\sF)$ reduces to
\[
0 \to F_{-l} \xrightarrow{\psi_{-l+1}} F_{-l+1} 
\xrightarrow{\psi_{-l+2}}
\dots \xrightarrow{\psi_3} F_3 \xrightarrow{\psi_4} F_4 \to 0.
\]
Notice furthermore that each cohomology group $H^i(\sI_S(k))$ 
appears in
at most one $F_e$.

In the following tables the entry $F_e$ in the row of $h^j$ and
column $k$ indicates that $H^j(\sI_S(k))\otimes \bigwedge^{n-k} 
\sU_l$ is a summand of $F_e$.
For convenience we indicate $\bigwedge^{n-k} \sU_l$ in the first row.

\[
\begin{matrix}
\text{For hyperplanes ($l=1$):} && \text{For planes ($l=2$):} \\
\\
\begin{array}{c|c|c|}
         & \sO(-1) & \sO  \\
\hline
     h^4 & F_3 & F_4    \\
\hline
     h^3 & F_2 & F_3    \\
\hline
   h^2  &F_1 & F_2    \\
\hline
 h^1 & F_0 & F_1    \\
\hline
h^0  & F_{-1} & F_0    \\
\hline
  & n-1 & n    \\
 \end{array}
&&
\begin{array}{c|c|c|c|}
        &  \sO(-1) & \sU_2 & \sO \\
\hline
h^4 &   F_2 & F_3 & F_4   \\
\hline
h^3 &   F_1 & F_2 & F_3    \\
\hline
h^2 &   F_0 & F_1 & F_2    \\
\hline
h^1 &   F_{-1} & F_0 & F_1    \\
\hline
h^0 &   F_{-2} & F_{-1} & F_0   \\
\hline
   &  n-2 & n-1  & n    \\
\end{array}
\\
\\
\\
\text{For lines ($l=3$):} && \text{For points ($l=4$):}\\
\\
\begin{array}{c|c|c|c|c|}
&   \sO(-1) & \bigwedge^2 \sU_3 & \sU_3 & \sO \\
\hline
h^4 &  F_1 & F_2 & F_3 & F_4    \\
\hline
h^3 &  F_0 & F_1 & F_2 & F_3    \\
\hline
h^2 &  F_{-1} & F_0 & F_1 & F_2    \\
\hline
h^1 &  F_{-2} & F_{-1} & F_0 & F_1    \\
\hline
h^0 &  F_{-3} & F_{-2} & F_{-1} & F_0    \\
\hline
 &  n-3 & n-2 & n-1 & n    \\
\end{array}
&&
\begin{array}{c|c|c|c|c|c|c|c}
&   \sO(-1) & \bigwedge^3 \sU_4 & \bigwedge^2 \sU_4 & \sU_4 & \sO \\
\hline
h^4 &  F_0 & F_1 & F_2 & F_3 & F_4    \\
\hline
h^3 &  F_{-1}& F_0 & F_1 & F_2 & F_3    \\
\hline
h^2 &  F_{-2} & F_{-1} & F_0 & F_1 & F_2    \\
\hline
h^1 &  F_{-3} & F_{-2} & F_{-1} & F_0 & F_1    \\
\hline
h^0 &  F_{-4} & F_{-3} & F_{-2} & F_{-1} & F_0    \\
\hline
 &  n-4 & n-3 & n-2 & n-1 & n    \\
\end{array}
\\
\end{matrix}
\]

We mainly use this setup to calculate the cohomology of hyperplane, 
plane and line sections of $S$ 
via the following propositions:

\begin{prop} \xlabel{pCohomology}
Let $\sigma \in \GG_l$ be a linear subspace of codimension $l$. 
If $\psi_{i+1}= 0$  then
\[
H^i(\sI_S(n)|_{\sigma}) = \coker \psi_{i}.
\]
In particular if  $F_i = 0$  then $H^i(\sI_S(n)|_{\sigma}) =0$ for all
$\sigma \in \GG_l$.
\end{prop}

\begin{proof}
Restriction is right exact. 
\end{proof}

\begin{rem}
Notice that this proposition gives in a compact way the information 
one would get by
repeatedly using the restriction sequence.
\end{rem}

We compare the restriction of the ideal sheaf to the ideal sheaf 
of the restriction:

\begin{prop} \xlabel{pExcess}
Let $\sigma \in \GG_l$ be a linear subspace of codimension $l$, which 
does not intersect $S$ properly, but is not contained in $S$. Let $Z 
\subset S \cap \sigma$
the union of those components that are of codimension $1$ in 
$\sigma$. Then there exists a locally free sheaf $J_{Z}$ on $Z$ such 
that 
\[
0 \to J_{Z} \to I_S|_\sigma \to I_{S\cap \sigma} \to 0.
\]
Furthermore $J_Z$ completes the diagram
$$
\begin{matrix} &&0&&0&&&&\\
&&\downarrow &&\downarrow&&&&\\
0&\to &J_{Z} &\to &N^*_{\sigma/\PP^4}|_{Z}&\to&N^*_{Z/S}&\to &0\\
&&\downarrow &&\downarrow&&\| &&\\
0&\to &N^*_{S/\PP^4}|_{Z} &\to &N^*_{Z/\PP^4}&\to&N^*_{Z/S}&\to&0\\
&&\downarrow &&\downarrow&&&&\\
0&\to &N^*_{Z/\sigma}&= &N^*_{Z/\sigma}&&&&\\
&&\downarrow &&\downarrow&&&&\\
&&0 &&0&&&&\\
\end{matrix}
$$
of conormal sheafs. 
\end{prop}

\begin{proof}
Locally at a point $z\in Z$, let $I=\sI_{S,z}$ be the ideal of $S$ 
and 
$J=\sI_{\sigma ,z}$ the ideal of $\sigma$. Then 
the 
restriction  $I|_\sigma$ is given by $\frac I{IJ}$, 
since it is the ideal $I$ tensored by the local coordinate ring 
$R_{\sigma}$ of 
$\sigma$ at $z$.   
The ideal of the intersection is $\frac {I+J}{J}$ inside 
$R_{\sigma}$.   Furthermore there is a 
natural 
surjective map
 $\frac I{IJ}\to \frac {I+J}{J}$ of $R_{\sigma}$ modules.  The kernel 
is 
easily identified in the exact sequence:
$$0\to \frac {I\cap J}{IJ}\to \frac I{IJ}\to \frac {I+J}{J}\to 0.$$
Notice that the kernel is supported where the intersection is not 
proper, i.e. on $Z$.
In particular there is a short exact sequence
$$0\to J_{Z}\to \sI_{S}|_\sigma\to \sI_{S\cap \sigma}\to 0.$$

Now since $Z$ has pure codimension $1$ in $\sigma$ it is a local complete intersection. Also 
$\sigma$ and $S$ are smooth, so we have an exact sequence
\[
0\to \frac {I\cap J}{IJ}\to \frac{I}{I^2}\otimes 
\frac{R_{\PP^4,z}}{J} 
\to \frac {J}{J^2}\to 0.
\]
This proves that $J_Z$ is locally free and fits into the proposed 
diagram. 
\end{proof}

\begin{rem}
In Fultons notation \cite[Section 6.3]{FultonIntersection}, the dual 
of $J_{Z}$ is called 
the excess normal bundle of the fiber product 
$$
\begin{matrix}  Z&\to &S\\
\downarrow &&\downarrow\\
\sigma &\to &\PP^4 &\!\!\!\!\!.\\
\end{matrix}
$$
\end{rem}

\section{Lines} \xlabel{sLines}
\nosubsections

We consider multi-secants and lines in $S$.

\begin{prop} \xlabel{pSecantH}
Let $L$ be a $k$-secant line to $S$. Then 
\[
h^1(\sI_S(n)|_L) = 
\left\{
\begin{matrix}
0 & \text{if $n>k-2$}\\
k-n-1 & \text{ if $n\le k-2$}
\end{matrix}
\right.
\]
\end{prop}

\begin{proof} 
We are in the situation of Propositon \ref{pExcess} with $Z$ a scheme 
of length $k$ and an exact sequence
\[
0 \to J_Z \to \sI_S(n)|_L  \to  \sI_Z(n)\to 0.
\]
Since $J_Z$ has no $H^1$ we obtain $h^1(\sI_S(n)|_L)= h^1(\sI_Z(n)) = 
h^1(\sO_{\PP^1}(n-k))$.
\end{proof}

\begin{prop} \xlabel{pMinusoneline}
Let $L\subset S$ be a $(-k)$-line. Then 
\[
\sI_S|_L  = \sO(-a) \oplus \sO(-b)
\]
with $a+b=k+3$ and $0 < a,b < k+3$.
\end{prop}

\begin{proof}  If $L$ is contained in $S$, then $\sI_S|_L $ is the 
conormal bundle of $S$ restricted to $L$ and fits into the conormal 
bundle sequence 
\[
0 \to N_{S/\PP^4}^*|_{L} \to N_{L/\PP^4}^* \to N^*_{L/S} \to 0
\]
which reduces to
\[
0 \to N_{S/\PP^4}^*|_{L}  \to 3 \sOPone(-1) \to \sOPone(k) \to 0,
\]

from which the proposition follows.

\end{proof}

\begin{cor} \xlabel{cLineInSurfaceH}
Let $L \subset S$ be a $(-k)$-line.
If $n \ge k+1$ then 
\[
h^0(\sI_S(n)|_L)-h^0(\sO_L(n))= n-k-2.
\] \xlabel{coNmanysections}
On the other hand, if $h^1(\sI_S(n)|_L) \not=0$ then $k \ge n$. 
\xlabel{coNfewsections}
\end{cor}

\begin{proof}
 $\sI_S(n)|_L = \sO(n-a)\oplus \sO(n-b)$ with
$a,b<k+3$ by Proposition \ref{pMinusoneline}. For $n \ge k+1$ we have 
$h^1(\sI_S(n)|_L)=0$ and
$h^0(\sI_S(n)|_L)=2n-k-1$ since $a+b=k+3$.
\end{proof}

\begin{prop} \xlabel{pRestrictToLine}
Let $L \subset \PP^4$ be any line. 
\begin{enumerate}
\item If $h^1(\sI_S(n)|_L) = 1$ then either $L$ is a $n+2$-secant 
line or $L \subset S$ and
$L^2 \le -n$. \xlabel{linewithh1}
\item If $h^1(\sI_S(n)|_L) = 0$ and $h^0(\sI_S(n)|_L) - 
h^0(\sO_L(n)) = 1$ then
$L \subset S$ with $L^2 = 3-n$. \xlabel{linewithouth1}
\end{enumerate}
\end{prop}

\begin{proof}
 From Proposition \ref{pSecantH}  and 
Corollary \ref{cLineInSurfaceH} we obtain (\ref{linewithh1}).
For $L \not\subset S$ one always 
has
$h^0(\sI_S(n)|_L) \le h^0(\sO_L(n))$, so claim (\ref{linewithouth1}) follows.
\end{proof}

\section{Planes} \xlabel{sPlanes}
\nosubsections

Throughout this section
let $P \subset \PP^4$ be a plane and $C \subset P \cap S$ be the 
$1$-dimensional component. By Proposition \ref{pExcess} we have
the exact sequence
$$0\to J_{C}\to N^{*}_{P/\PP^4}|_C\to N^*_{C/S}\to 0$$
which reduces to
$$0\to J_{C}\to 2\sO_{C}(-1)\to \sO_{C}(-C)\to 0$$
and shows $J_{C}=\sO_{C}(C-2H)$. 

We can read off the degree of $C$ from the cohomology of either 
$I_{S\cap P}$ or $I_S|_P$:

\begin{prop}\xlabel{pPlaneH2} Let $P \subset \PP^4$ be any plane and
$C$ the curve component of $P \cap S$.
\begin{enumerate}
    \item If $n \ge -2$ then
$
h^2(\sI_{S}(n)|_P)=h^2(\sI_{P \cap S}(n)) =  { {\deg C - n - 1} 
\choose 
{2}}.
$

\item If $h^1(\sO_{C}(C+(n-2)H))=0$ then 
$
h^1(\sI_{S}(n)|_P)=h^1(\sI_{P \cap S}(n)).
$
\end{enumerate}
\end{prop}

\begin{proof}
First we consider the cohomology of the short exact sequence
$$0\to \sO_{C}(C+(n-2)H)\to \sI_{S}(n)|_P\to \sI_{S\cap P}(n)\to 0.$$
Since $h^2(\sO_{C}(C+(n-2)H))=0$, we get 
$h^2(\sI_{S}(n)|_P)=h^2(\sI_{P \cap S}(n))$.  In the second case
$h^1(\sO_{C}(C+(n-2)H))=h^2(\sO_{C}(C+(n-2)H))=0$ and the second part 
of the proposition follows.

From the sequence
\[
0 \to \sI_{P \cap S}(n) \to \sO_P(n) \to \sO_{P \cap S}(n) \to 0
\]
we obtain $h^2(\sI_{P \cap S})(n) = h^1(\sO_{P \cap S}(n))$ if $n \ge 
-2$. Furthermore
\[
h^1(\sO_{P \cap S}(n)) = h^1(\sO_{C}(n)) = h^0(\sO_C(\deg C-3-n)) 
=  {\deg C - n  - 1\choose 2}
\]

\end{proof}

\section{diagrams} \xlabel{sDiagrams}
\nosubsections
In the previous section we compared the restriction of the ideal of a 
surface in $\PP^4$ to a plane with the ideal of the intersection of 
the surface and the plane.  In this section we will concentrate on 
the latter. We recall from \cite{GreenGin} how the different Hilbert 
functions of plane algebraic sets are read off from certain diagrams 
representing the generic initial ideals of their ideals.

\begin{notation} Consider

\begin{tabular}{ll}
$K[a,b,c]$ & the coordinate ring of $\PP^2$ \\
 $\gin I$ & the generic initial ideal of $I$ with respect to the \\ & 
reverse lexicographic order
                   with $a<b<c$.
\end{tabular}
\end{notation}

\begin{rem} \quad

\begin{enumerate}
\item $I$ is saturated if and only if $\gin I$ is saturated. 
\item if $\gin I$ is saturated an $a^ib^jc^k \in \gin I$ then also 
$a^ib^j \in \gin I$.
\item The Hilbert function and Hilbert polynomial of $I$ and $\gin I$ 
are the same.
\end{enumerate}
\end{rem}

\begin{defn}
Let $I \subset K[a,b,c]$ be a saturated ideal. We represent the 
generic initial ideal $\gin I$ by 
a diagram of $x$'s and $0$'s in $\NZ_{0}\times \NZ_{0}$. An $x$ in 
the 
point $(i,j)$ means 
that               
$a^ib^j \in \gin I$ and a $0$ means $a^ib^j \not\in \gin I$.

We also set 
\begin{align*}
d(I) &= \min\{i \suchthat \exists a^ib^j \in \gin I\} \\
e(I) &= \#\{a^ib^j \not\in \gin I \suchthat i\ge d\}
\end{align*}

\end{defn}

\begin{example}
$\gin I = (a^4,a^3b,a^2b^3)$ is represented by
\[
\begin{array}{c|cc|ccccccc}
& \vdots &\vdots & & \vdots \\
& 0 & 0 & x & x & x \\
& 0 & 0 & x & x & x  \\
j & 0 & 0 & 0 & x & x & \hdots\\
& 0 & 0 & 0 & x & x  \\
& 0 & 0 & 0 & 0 & x  \\
\hline
& & & & i & &\\
\end{array}
\]
We have $d(I)=2$ and $e(I)=4$. Notice that $d(I)$ is the number of 
columns with only $0$'s
and $e$ is the number of $0$'s outside of these columns.
\end{example}

\begin{rem}
Since $\gin I$ is an ideal we have $0$'s left and below each 
$0$, and $x$'s
right and above of each $x$.
Also generic initial ideals are Borel-fixed, i.e. for $i \ge 1$ we 
have
\[
a^ib^j \in \gin I \implies a^{i+1}b^{i-1} \in \gin I.
\]
This means that we also have $x$'s on the diagonal right and below
of each $x$.
\end{rem}

\begin{prop} \xlabel{pBelow}
Let $\sI$ be an ideal sheaf on $\PP^2$ and $I = \bigoplus H^0(\sI(n))$
the corresponding saturated ideal. Let  
$H_{I}(n)=h^0(\sO_{\PP^2}(n))-h^0(\sI(n))$ be the  Hilbert function
of $V(I)$. Then
\[
H_{I}(n) = \# \{ \text{$0$'s below and on
the diagonal $i+j=n$ in the diagram of $\gin I$}\}.
\]
\end{prop}

\begin{proof}
The monomials not 
contained
in $\gin I$ form a basis of $K[a,b,c]/I$. 
\end{proof}

\begin{prop}
Let $I \subset K[a,b,c]$ be a saturated ideal. Then 
$d(I)$ as defined above is the degree of the curve components of 
$V(I)$.
\end{prop}

\begin{proof}
For large $n$ the number of $0$'s on and below $i+j=n$ increases by 
$d(I)$ in each step, so the linear term of the Hilbert polynomial  of 
$V(I)$ has 
coefficient $d(I)$.
\end{proof}

\begin{prop}
Let $I \subset K[a,b,c]$ be a saturated ideal. Then 
$e(I)$ as defined above is the degree of the dimension $0$ component 
of $V(I)$.
\end{prop}

\begin{proof} The difference of the constant in the Hilbert 
polynomial 
of $V(I)$ and that of a plane curve of degree $d(I)$ is precisely $e$.
\end{proof}

Useful for the geometric interpretation is the following

\begin{rem} 
It follows from a theorem of Ellia and Peskine \cite[Theorem 
4.4]{GreenGin} that 
one 
can sometimes
read off special positions of points in $V(I)$ according to the 
following rule:

If in the $i$'th column we have at least three $x$'s to the right of 
three $0$'s in column $i-1$,
then there exists a curve of degree $i-d(I)$ passing through $n$ 
points of $V(I)$, where 
$n$ is the number of $0$'s in columns $d(I)+1,\dots,i-1$.  The 
converse of this is
not true in general.
\end{rem}

\begin{prop}
Let $\sI$ be an ideal sheaf on $\PP^2$, $I = \bigoplus H^0(\sI(n))$
the corresponding saturated ideal and $r$ the number of $0$'s
with $i\ge d(I)$ lying above the diagonal $i+j=n$,
then
\[
h^1(\sI(n)) = r
\]
\end{prop}

\begin{proof} The number $r$ is the difference between the Hilbert 
    function and the Hilbert polynomial of $V(I)$ at $n$. 
\end{proof}

We again turn to our smooth surface $S$ of degree 11 in $\PP^4$.
\begin{rem}
Since the degree of the curve component of a plane section $S \cap P$
is bounded by the degree of $S$,
it follows that for $h^1(\sI_{S\cap P}(n))=r$ there are only 
finitely many saturated generic initial ideals.
\end{rem}

\begin{example}
Let $X\subset \PP^2$ be a finite subscheme of degree $11$ that is not 
contained in any conic section.
For the diagram of $X$ this means
that we must have no $x$'s on and below the $i+j=2$ line, and eleven 
$0$'s altogether.  The possible such diagrams
are
\[
\begin{array}{|ccccccccc}
0 & x & x & x   \\
0 & x & x & x   \\
0 & x & x & x   \\
0 & x & x & x   \\
0 & x & x & x   \\
0 & x & x & x \\
0 & 0 & x & x  \\
0 & 0 & 0 & x    \\
\hline
\end{array} \quad \quad
\begin{array}{|ccccccccc}
0 & x & x & x   \\
0 & x & x & x   \\
0 & x & x & x   \\
0 & x & x & x \\
0 & 0 & x & x \\
0 & 0 & x & x  \\
0 & 0 & 0 & x    \\
\hline
\end{array}  \quad \quad
\begin{array}{|ccccccccc}
0 & x & x & x \\
0 & x & x & x \\
0 & 0 & x & x \\
0 & 0 & x & x \\
0 & 0 & x & x  \\
0 & 0 & 0 & x    \\
\hline
\end{array}
\]
and
\[
\begin{array}{|ccccccccc}
0 & x & x & x   \\
0 & x & x & x   \\
0 & x & x & x   \\
0 & 0 & x & x \\
0 & 0 & 0 & x  \\
0 & 0 & 0 & x    \\
\hline
\end{array} \quad \quad
\begin{array}{|ccccccccc}
x & x & x & x   \\
0 & x & x & x   \\
0 & 0 & x & x \\
0 & 0 & x & x \\
0 & 0 & 0 & x  \\
0 & 0 & 0 & x    \\
\hline
\end{array}  \quad \quad
\begin{array}{|ccccccccc}
x & x & x & x \\
0 & x & x & x \\
0 & x & x & x \\
0 & 0 & x & x \\
0 & 0 & 0 & x  \\
0 & 0 & 0 & 0    \\
\hline
\end{array}
\]
In the first case $X$ contains a subscheme of length $8$ on a line, 
in the second case $X$ contains a subscheme of length $7$ on a line, 
in the third case $X$ contains a subscheme of length $10$ on a conic, 
in the fourth case $X$ contains a subscheme of length $6$ on a line, 
in the first five cases $X$ is contained in a cubic, while in the 
last case $X$ is not contained in any cubic curve.  Notice that each 
case is also distinguished by the corresponding values of 
$h^1(I_X(n))$ for $n=3,4,5$. In fact, we get the following triples 
$(4,3,2)$, $(3,2,1)$, $(3,1,0)$, $(2,1,0)$, $(2,0,0)$ and $(1,0,0)$ 
respectively.
\end{example}

\begin{example} \xlabel{eQuartic3}
Let $X\subset \PP^2$ be the union of a quartic curve and a scheme of 
length $3$.  In the diagram of $X$ this means that the first four 
columns have all $0$'s, and that there are three more $0$'s.  There 
are two possible diagrams:

\[
\begin{array}{|cccc|ccccc}
0 & 0 & 0 & 0 & x & x\\
0 & 0 & 0 & 0 & x & x\\
0 & 0 & 0 & 0 & x & x\\
0 & 0 & 0 & 0 & x & x \\
0 & 0 & 0 & 0 & 0 & x\\
0 & 0 & 0 & 0 & 0 & x \\
0 & 0 & 0 & 0 & 0 & x  \\
\hline
\end{array} \quad
\begin{array}{|cccc|ccccc}
0 & 0 & 0 & 0 & x & x\\
0 & 0 & 0 & 0 & x & x\\
0 & 0 & 0 & 0 & x & x\\
0 & 0 & 0 & 0 & x & x \\
0 & 0 & 0 & 0& x & x\\
0 & 0 & 0 & 0 & 0 & x  \\
0 & 0 & 0 & 0 & 0 & 0  \\
\hline
\end{array}
\]
In both cases we have
$h^1(I_X(4))=3$.  But only in the first case is the length 
$3$-subscheme on a line and only in this case is $h^1(I_X(6))=1$.
\end{example}

\section{Rational Surfaces with $d=11$, $\pi=11$} \xlabel{sElevenEleven}
\nosubsections

Let $S \subset \PP^4$ be a rational surface of degree $d=11$ an
section genus $\pi=11$. In this section we determine the possible 
Hilbert functions
that $S$ can have and find restrictions on the maps in the Tate 
resolution.
By Popescu \cite{PopescuThesis} such surfaces have the following 
cohomology table for the ideal sheaf $\sI_S$

\[
\begin{array}{c|c|c|c|c|c|c|c|c|}
\hline
h^4 & & &    &    &    &          & &\\ \hline
h^3 & 11 & &    &    &    &          && \\ \hline
h^2 && & 3 & 1 &    &         &&\\ \hline
h^1 &&  &   &    & 2 & 1+a & b & c \\ \hline
h^0 &&  &   &    &    & a      & 10+b & 38+c\\ \hline
&  \sI_S(-1) & \sI_S(0) &  \sI_S(1) & \sI_S(2) & \sI_S(3)  & \sI_S(4) 
& \sI_S(5)  & \sI_S(6)\\
\end{array}
\]
We now consider the Tate resolution of $\sI_{S}(n)$. The most 
interesting part
for our purposes is
\[
\cdots
\longto{}
3E(n-1)
\longto{}
\begin{matrix}
E(n-2) \\ \oplus \\
2E(n-3) \\ \oplus \\
aE(n-4)
\end{matrix}
\longto{}
\begin{matrix}
(a+1) E(n-4) \\ \oplus \\
(10+b) E(n-5)
\end{matrix}
\longto{}
\begin{matrix}
b E(n-5) \\ \oplus \\
(38+c) E(n-6)
\end{matrix}
\longto{}
\cdots
\]
and in particular the subcomplex
\[
3E(n-1)
\longto{\alpha}
\begin{matrix}
E(n-2) \\ \oplus \\
2E(n-3)
\end{matrix}
\longto{\beta}
\begin{matrix}
(a+1) E(n-4)
\end{matrix}
\longto{\gamma}
\begin{matrix}
b E(n-5)
\end{matrix}.
\]
Applying the functor $U_l$ for $l=1,2,3$ we get complexes 
$U_l(\sI_S(n))$ with maps
 $U_l(\alpha)$, $U_l(\beta)$ and $U_l(\gamma)$.  By abuse of notation 
we often drop the functor. First we use Proposition
\ref{pCohomology} to compute the possible cohomology groups for the 
restriction 
of $\sI(n)$ to linear subspaces:

\begin{prop} \xlabel{pCohomologyEleven}
Let $\sigma \in \GG_l$ be a linear subspace of codimension $l=1,2,3$.
Then the cohomology
table of $\sI_S(n)|_{\sigma}$ for $n=2,3,4$ may have the following 
entries:
\[
\begin{matrix}
\text{for $l=1$ (hyperplanes):} & \text{for $l=2$ (planes):} &
\text{for $l=3$ (lines):} \\ \\
\begin{array}{|c|c|c|c|c|c|c|c|c|}
\hline
&    &       \\ \hline
&    &       \\ \hline
0/1&     &       \\ \hline
2/3 & 3 & a-1/a/a+1     \\ \hline
      &    &  a/a+1/a+2   \\ \hline
\end{array}
&
\begin{array}{|c|c|c|c|c|c|c|c|c|}
\hline
&    &       \\ \hline
&    &       \\ \hline
0/1&     &       \\ \hline
? & 1\dots4 & a-4 \dots a+1     \\ \hline
   &  ?  & ?   \\ \hline
\end{array}
&
\begin{array}{|c|c|c|c|c|c|c|c|c|}
\hline
&    &       \\ \hline
&    &       \\ \hline
0 &     &       \\ \hline
? & 0\dots5 & a-8 \dots a+1     \\ \hline
   &  ?  & ?   \\ \hline
\end{array}
\\
\end{matrix}
\]
where empty boxes stand for cohomologies that must be zero and 
question marks stand for cohomologies
for which we have no restrictions so far.
\end{prop}

\begin{proof}
We use the first part of Proposition \ref{pCohomology} repeatedly,
and indicate the ranks of the vector bundles in the source and in the 
target
to find the possible coranks of the maps.
The cohomology group $H^2(\sI_S(2)|_\sigma)$ is the cokernel of 
$\fU_l(\alpha)$
whose source and target have ranks $3 \to 1$, $3\cdot 2 \to 1$ and $3 \cdot 3 \to 1$ for
$l=1,2,3$ respectively. In addition we must have
$h^2(\sI_S(2)|_\sigma)=0$ for lines.

Similarily the cohomology groups $H^1(\sI_S(3)|_\sigma)$ are 
cokernels of $\fU_l(\alpha)$ whose source and target have ranks $0 \to 1+2$, $3 \to 1\cdot 2+2$ 
and $3\cdot 3 \to 1\cdot 3 + 2$ respectively. Since the map $a \sO \to (a+1) \sO$ 
is always zero, the 
cohomology groups  $H^1(\sI_S(4)|_\sigma)$ are cokernels of
$\fU_l(\beta)$ whose source and target have ranks $2 \to
1+a$, $1+2\cdot2 \to 1+a$ and $1\cdot 3 + 2 \cdot 3 \to 1+a$ 
respectively.

For hyperplanes the intersection $S\cap \sigma$ is always a curve of degree $11$ and 
arithmetic genus $11$. The possible cohomology dimensions  
$h^1(\sI_S(2)|_\sigma)$, $h^0(\sI_S(3)|_\sigma)$ and 
$h^0(\sI_S(4)|_\sigma)$ are therefore determined by Riemann-Roch.

The empty boxes of the proposition follow from the second part of 
Proposition \ref{pCohomology}.
\end{proof}

Consider the linear part $\alpha_{1}$ of $\alpha$ in the Tate 
resolution.   It is
given by  a
$(3\times 1)$ matrix with entries in
$W^{*}$. These entries can be interpreted as points in $\PP^4$.

\begin{prop}\xlabel{plane5}
Let $\sigma$ be the linear space spanned by the entries of 
$\alpha_{1}$ in the Tate resolution of $\sI$.
 Then
$\sigma=P$ is a plane and $P\cap S$ contains the
unique plane quintic curve on $S$.
\end{prop}

\begin{proof}
If $\sigma=P$ is a plane, we consider the
map
\[
3 \sU_2 \xrightarrow{\alpha_{1}} \sO
\]
on $\GG_2$. It drops rank only on $P
\in \GG_2$. By Proposition
\ref{pCohomology} and $\ref{pPlaneH2}$ this happens if and only
if $P\cap S$ contains a plane quintic.

If $\sigma$ is not a plane, we choose a line $L$ that contains 
$\sigma$ and 
consider the map
\[
3 \sU_3 \xrightarrow{\alpha_{1}} \sO
\]
on $\GG_3$. If we restrict to $L$, this map vanishes and 
we obtain $h^2(\sI_S(2)|_\sigma)=1$ by Proposition 
\ref{pCohomology}. This is impossible on a line. 
\end{proof}

Let $C$ be the unique plane quintic curve on $S$, let $P$ be its 
span, and let $D=H-C$ be the residual curve to $C$ in a hyperplane 
section. 
 Then $|D|$ is a pencil and $D^2$ is the length of the subscheme $R$ 
residual to $C$ in
 $P \cap S$.  

 \begin{lem} \xlabel{Dsquare}
$0\leq D^2\leq 2$ and a general member of $|D|$ is a smooth curve of 
genus $g(D)=D^2$.
\end{lem}
\begin{proof}  Since $|D|$ has no fixed component $D^2\geq 0$.  It 
remains to show that $D^2\leq 2$.

Consider the short exact sequences of Proposition \ref{pExcess}:
$$0\to J_{C}(3)\to \sI_{S}|_P(3)\to \sI_{S\cap P}(3)\to 0.$$
and 
$$0\to J_{C}(3)\to 2\sO_{C}(2)\to \sO_{C}(3H-C)\to 0.$$

Notice that $h^0(\sI_{S\cap P}(3))=0$, 
so taking cohomology in the former sequence yields  
$h^1(\sI_{S}|_P(3))=h^1( J_{C}(3))+h^1(\sI_{S\cap P}(3))$.  
Furthermore $h^1(\sI_{S\cap P}(3))=h^1(\sI_{R}(-2))=D^2$.  Therefore 
 $$h^1(\sI_{S}|_P(3))= h^1( J_{C}(3))+h^1(\sI_{S\cap P}(3))=h^1( 
J_{C}(3))+D^2.$$
 On the other hand $h^1(\sI_{S}|_P(3))\leq 4$ by Proposition 
\ref{pCohomologyEleven}, so $h^1( J_{C}(3))+D^2\leq 4$.
First, this implies that $D^2\leq 4$, which means that 
$D \cdot C =(H-D)D=6-D^2 \ge 2$.  
But $\sO_{C}(3H-C)=\sO_C(2H+D)=\omega_C(D)$, so 
$h^1(\sO_{C}(3H-C))=0$.
Secondly, taking cohomology in the second sequence we get $h^1( 
J_{C}(3))\geq 
2h^1(\sO_{C}(2))=2$, so $2+D^2\leq h^1( J_{C}(3))+D^2\leq 4$, i.e. 
$D^2\leq 2$.

The pencil of curves $|D|$ has a base locus of length at most $2$. By 
Bertini's Theorem the general member has singularities only in this
base locus.  
But if  the general $D$  is singular in the  base locus, then 
$D^2\geq 4$, so we conclude that $D$ is smooth.
Furthermore, $|D|$ is complete as a linear system, in fact $|H|$ is 
complete by Severis Theorem and embeds $C$, so  $D=H-C$ can only move 
in a pencil.  But the general member of a complete pencil of curves 
on a rational surface, that does not have a fixed component, must be 
irreducible:  In fact, the connected fibers of the Stein 
factorization of the map defined by $|D|$ are already linearly 
equivalent.  Consequently, if the general element $D$ is a multiple 
of fibers, one could move one fiber while fixing the rest, 
contradicting the assumption that $|D|$ has no fixed component.
Therefore the general member $D$ of the pencil $|D|$ is a
smooth and irreducible curve of genus 
$$g(D)=\frac 12 (D^2+D\cdot K)+1=\frac 12 (D^2-2+D^2)+1=D^2.$$
\end{proof}

We take a closer look at the subcomplex
\[
3E(n-1)
\longto{\alpha}
\begin{matrix}
E(n-2) \\ \oplus \\
2E(n-3)
\end{matrix}
\longto{\beta}
\begin{matrix}
(a+1) E(n-4)
\end{matrix}
\longto{\gamma}
\begin{matrix}
b E(n-5)
\end{matrix}.
\]
of the Tate resolution of $\sI_{S}(n)$.
The maps $\alpha$, $\beta$ and $\gamma$ can be given by
matrices with entries of the following degrees:
\[
A = \begin{pmatrix} 1 & 2 &2 \\ 1 & 2 & 2  \\ 1 & 2 & 2 \end{pmatrix}
\quad 
B = \begin{pmatrix}
2 & \cdots & 2 \\
1 & \cdots & 1 \\
1 & \cdots & 1 \\
\end{pmatrix}
\quad
\Gamma = \begin{pmatrix}
1 & \cdots & 1 \\
\vdots && \vdots \\
1 & \cdots & 1 \\
\end{pmatrix}.
\]
Notice that $A$, $B$ and $\Gamma$ do not depend on the twist $n$.

\begin{prop} \xlabel{pColumnsAlgebraic}
Let $b = (q,l_1,l_2)^T$ be an column vector over the exterior algebra
$E = \bigwedge V$ with entries of degress $(2,1,1)^T$. Then after
coordinate changes and row operations we have one of following
possiblities:
\begin{enumerate}
\item $b = (e_3\wedge e_2,e_1,e_0)^T$ \xlabel{iGeneral}
\item $b = (0,e_1,e_0)^T$ \xlabel{iOneSixsecant}
\item $b = (e_4\wedge e_3+e_2 \wedge e_1,e_0,0)^T$ \xlabel{iFullRank}
\item $b = (e_2 \wedge e_1,e_0,0)^T$ \xlabel{iInfSixsecants}
\item $b = (0,e_0,0)^T$ \xlabel{iLinear}
\item $b = (e_3\wedge e_2 + e_1 \wedge e_0,0,0)^T$ \xlabel{iQ2}
\item $b = (e_1 \wedge e_0,0,0)^T$ \xlabel{iQ1}
\item $b = (0,0,0)^T$ \xlabel{iQ0}
\end{enumerate}
with $e_0\dots e_4$ a basis of $ V = W^*$.
\end{prop}

\begin{proof}
We collect the coefficients of $q$ in a skew symmetric $5 \times 5$ 
matrix $M$. We say that $q$ has rank $r$ if $M$ has rank $2r$. 

If the linear forms are independent, we can assume that $q$ involves 
only the remaining $3$ variables. Consequently we have $\rank q \le 
1$. This gives the cases (\ref{iGeneral}) and (\ref{iOneSixsecant}).

If the linear forms span a $1$ dimensional space, we can assume that 
$q$ involves only the remaining $4$ variables and $\rank q \le 2$. 
This gives the cases (\ref{iFullRank}), (\ref{iInfSixsecants}) and 
(\ref{iLinear}).

If both linear forms are zero, $q$ can involve all $5$ variables and 
$\rank q \le 2$. This gives
the last three cases
\end{proof}

\begin{prop} \xlabel{pSecants}
Let $b$ be a column of $B$. Then after row operations and coordinate 
changes one
of the following holds
\begin{enumerate}
\item $b = (e_3\wedge e_2, e_1,e_0)^T$ and  $P$ is contained in the 
$\PP^3$ spanned by
 $e_0$, $e_1$, $e_2$ and $e_3$. Furthermore the line $L$ though $e_0$ 
and $e_1$ either intersects $S$ in a scheme of length at least $5$ or 
$L \subset S$ 
with $L^2 \le -3$. If $L$ lies in $P$ then $b$ vanishes on $P$. 
\xlabel{iGeneralGeometry}
\item $b = (e_2\wedge e_1, e_0,0)^T$ and $P$ is spanned by $e_0$, 
$e_1$ and $e_2$. Furthermore each line $L$ that passes 
through $e_0$ and lies in $P$  either intersects $S$ in a scheme of length at least $6$ or 
$L\subset S$ with $L^2 \le -4$.  
\xlabel{iSixInfGeometry}
\item $b = (0,e_1,e_0)^T$ and
the line $L$ though $e_0$ and $e_1$ either intersects $S$ in a scheme 
of length at least $6$
or $L \subset S$ with $L^2 \le -4$. \xlabel{iSixOneGeometry}
\end{enumerate}
\end{prop}

\begin{proof}
First we look at the three possible cases and afterwards we exclude 
all other possibilities in Proposition \ref{pColumnsAlgebraic}. In 
both parts we use on the one hand the fact that $AB =0$ in the 
exterior algebra to obtain information about the matrix $A$ and on 
the other hand the geometric interpretation of $A$ and $B$ for 
various $l$.

The syzygy matrix of  $b = (e_3\wedge e_2, e_1,e_0)^T$ is
$$\begin{pmatrix}0&
      0&
      0&
      {e}_{{3}}&
      {e}_{{2}}&
      {-{e}_{1}}&
      {-{e}_{0}}\\
      0&
      {e}_{1}&
      {e}_{0}&
      0&
      0&
      {e}_{{3}} \wedge {e}_{{2}}&
      0\\
      {e}_{0}&
      0&
      {e}_{1}&
      0&
      0&
      0&
      {e}_{{3}} \wedge {e}_{{2}}\\
      \end{pmatrix}^T$$
therefore the linear part of $A$ contains linear combinations of $e_0 
\dots e_3$. This proves that
$P$ lies in the $\PP^3$ spanned by these points in $\PP^4$. Since the 
line $L$ also lies in this $\PP^3$ it is either contained in $P$ and 
$A$ drops rank on $L$ or it intersects $P$ in a point $\lambda e_1 + 
\mu e_0$. On $L$ we then obtain
\[
A|_L =
\begin{pmatrix}
      {e}_{{3}}\wedge e_1 \wedge e_0&
      {e}_{{2}} \wedge e_1 \wedge e_0&
      0\\
      0&
      0&
      \lambda {e}_{{3}} \wedge {e}_{{2}} \wedge e_1 \wedge e_0\\
      0&
      0&
      \mu {e}_{{3}} \wedge {e}_{{2}} \wedge e_1 \wedge e_0\\
      \end{pmatrix}^T
\]
which also has submaximal rank. 
This implies the geometric properties of (\ref{iGeneralGeometry}) 
by Proposition \ref{pRestrictToLine}.
If $L$ lies in $P$ then $P= e_0 \wedge e_1\wedge (\lambda e_2 + \mu 
e_3)$ which annihilates all entries of $b$.  

The syzygy matrix of $b = (e_2\wedge e_1, e_0,0)^T$ is
$$\begin{pmatrix}0&
      0&
      {e}_{1}&
      {e}_{{2}}&
      {-{e}_{0}}\\
      0&
      {e}_{0}&
      0&
      0&
      {e}_{{2}} \wedge {e}_{1}\\
      1&
      0&
      0&
      0&
      0\\
      \end{pmatrix}^T$$
and $P$ is therefore spanned by $e_2$, $e_1$ and $e_0$. On the other 
hand $b$ vanishes on
all lines 
$L = (\lambda e_1+ \mu e_2) \wedge e_0$ and consequently $B$ drops 
rank there. This implies (\ref{iSixInfGeometry}) by Proposition 
\ref{pRestrictToLine}.

If $b = (0,e_1,e_0)^T$ this column vanishes on $L=e_1 \wedge e_0$ and 
this implies (\ref{iSixOneGeometry}).

Now we consider the other cases in Proposition 
\ref{pColumnsAlgebraic}.

The syzygy matrix  $b = (e_4\wedge e_3+e_2 \wedge e_1,e_0,0)^T$ is
$$\begin{pmatrix}0&
      0&
      {-{e}_{0}}&
      {e}_{{3}} \wedge {e}_{1}&
      {e}_{{4}} \wedge {e}_{1}&
      {e}_{{3}} \wedge {e}_{{2}}&
      {e}_{{4}} \wedge {e}_{{2}}&
      {e}_{{4}} \wedge {e}_{{3}}-{e}_{{2}} \wedge {e}_{1}\\
      0&
      {e}_{0}&
      {e}_{{4}} \wedge {e}_{{3}}+{e}_{{2}} \wedge {e}_{1}&
      0&
      0&
      0&
      0&
      0\\
      1&
      0&
      0&
      0&
      0&
      0&
      0&
      0\\
      \end{pmatrix}^T$$
and therefore $P$ must be spanned by $e_0$ which is impossible.

The column $b = (0,e_0,0)^T$ and $b=(0,0,0)$ vanish on $e_0$ which is 
impossible since no sheaf on 
a point can have nonvanishing $H^1$.

The column $b = (e_3\wedge e_2 + e_1 \wedge e_0,0,0)^T$ has syzygy 
matrix
$$\begin{pmatrix}0&
      0&
      {e}_{{2}} \wedge {e}_{0}&
      {e}_{{3}} \wedge {e}_{0}&
      {e}_{{2}} \wedge {e}_{1}&
      {e}_{{3}} \wedge {e}_{1}&
      {e}_{{3}} \wedge {e}_{{2}}-{e}_{1} \wedge {e}_{0}\\
      0&
      1&
      0&
      0&
      0&
      0&
      0\\
      1&
      0&
      0&
      0&
      0&
      0&
      0\\
      \end{pmatrix}^T$$
and there are no linear forms to span $P$.
      
Finally the syzygy matrix of $b = (e_1 \wedge e_0,0,0)^T$ is
$$\begin{pmatrix}0&
      0&
      {e}_{0}&
      {e}_{1}\\
      0&
      1&
      0&
      0\\
      1&
      0&
      0&
      0\\
      \end{pmatrix}^T$$
and $P$ must be spanned by $e_0$ and $e_1$ which is again impossible.
\end{proof}

\begin{rem}\xlabel{rsix}
Since $B$ has corank $1$ on any $6$-secant there is a unique column 
of type (\ref{iSixInfGeometry}) or (\ref{iSixOneGeometry}) for each 
such line. Furthermore $B$ vanishes on a $7$-secant line.
\end{rem}

\begin{rem}
Except for the observation that $A$ drops rank on $L$ in the first 
case, this classification was already obtained by Popescu in 
\cite{PopescuThesis}.
\end{rem}

\newcommand{\Pspan}{\PP_{\text{span}}}
\newcommand{\Pperp}{\PP_{\perp}}
\newcommand{\Pcol}{\PP^{a}_{c}}
\newcommand{\Pcolone}{\PP^{1}_{c}}
\newcommand{\Pcoltwo}{\PP^{2}_{c}}
\newcommand{\Prow}{\PP^1_{r}}
\newcommand{\pispan}{p}

We now consider the case of several columns in the matrix $B$, and 
start by focusing on the linear part $B_{1}$
and the span of its entries $\Pspan \subset \PP^4$.  We denote
by $\Pcol$ the column space and by $\Prow$ the row space of $B_1$. 
The Segre variety $\Prow \times \Pcol \subset \PP^{2a+1}$ is 
described by a $2\times (a+1)$ matrix and $B_1$ defines a birational 
map
\[
\pispan \colon \Prow \times \Pcol --\!\!> \Pspan\subset \PP^4
\]
which can be interpreted as the projection from a linear space 
$\Pperp \subset \PP^{2a+1}$ which is the space of linear relations 
between the entries of $B_1$. Denote by $T$ the image of $\pispan$.

\begin{lem} \xlabel{lContainingT}
Any quartic $X_4$ containing $S$ also contains $T$.
\end{lem}

\begin{proof}
If a column of $B$ is of type (\ref{iGeneralGeometry}) or 
(\ref{iSixOneGeometry}) then the span $L$ of its linear entries
is either contained in $S$ or at least a 
$5$-secant to $S$ by Proposition \ref{pSecants}. It is therefore 
contained in $X_4$. If a column of $B$ is of type 
(\ref{iSixInfGeometry}) it has only one linear entry which represents 
a point of 
$S \subset X_4$ by Proposition \ref{pSecants}.
\end{proof}

\begin{cor}\xlabel{Tlt2}
$\dim T \le 2$.
\end{cor}

\begin{proof}
If $\dim T > 2$ then by construction $a\ge2$ and $S$ is contained in 
at least $2$ independent quartics.
Since these also contain $T$ by Lemma \ref{lContainingT} we 
obtain a contradiction.
\end{proof}   

\begin{lem} \xlabel{lZfinite}
The intersection $Z = \Pperp \cap \Prow \times \Pcol$ is finite.
\end{lem}

\begin{proof}
A point $(r,c) \subset \Prow \times \Pcol$ is in $\Pperp$ if and only 
if the entry of $B_1$ in the corresponding generalized row and column 
is zero. If $Z$ is infinite, one of the following happens
\begin{enumerate}
\item a column of $B_1$ vanishes. This is impossible by Proposition 
\ref{pSecants}
\item several columns of $B_1$ have rank $1$. Each of them gives a 
residual point in $P$. 
\begin{enumerate}
\item If this point moves, we obtain infinitely many residual points 
in $P$ which is impossible.
\item If this point does not move, we obtain at least two columns 
that span only a point. Denote by $B'$ the corresponding two columns 
of $B$ and by $B_1'$ their linear part. Since there can be no zero 
column in $B_1'$ we obtain 
\[
B_1' = \begin{pmatrix} e_1 & 0 \\ 0 & e_1 \end{pmatrix}
\] 
after row and column operations.
By Proposition \ref{pSecants} the point $e_1$ lies in $P$. If $P$ is 
spanned by $e_1$, $e_2$ and $e_3$ the same proposition shows that the 
$2$-forms of $B'$ are linear combinations of $e_2 \wedge e_3$ and 
terms of the form $e_1 \wedge *$. After eliminating the $e_1$-terms 
with column operations we obtain
\[
B' = \begin{pmatrix}
\lambda e_2 \wedge e_3 & e_1 & 0 \\
\mu e_2 \wedge e_3 & 0 & e_1
\end{pmatrix}^T.
\]
By Proposition \ref{pSecants} again the coefficients $\lambda$ and 
$\mu$ 
must be nonzero, but then
we obtain the column $(0,\mu e_1, -\lambda e_1)^T$ as a linear 
combination contradicting Proposition \ref{pSecants}.
\end{enumerate}
\end{enumerate}
\end{proof}

\begin{cor} \xlabel{cPspanMin} In the above notation:

\begin{enumerate}
    \item
    $a\le \dim \Pspan$
    \item $a\le \dim T\le a+1$ \xlabel{iDimT}
    \item  $a\le 2$.
    \end{enumerate}
\end{cor}

\begin{proof}
If $\dim \Pspan < a$ we have $\codim \Pperp \le a$ and therefore 
$$\dim Z = 
\dim (\Pperp \cap \Prow \times \Pcol)\ge 1.$$
This 
contradicts Lemma \ref{lZfinite}.
The variety $T$ is the image of the projection from $\Pperp$ and 
$\Pperp \cap 
\Prow \times \Pcol$ is finite, so the fibers of the projection are at 
most $1$-dimensional and the second part follows.  The third part now 
follows from the second and Corollary \ref{Tlt2}.
\end{proof}

We now turn to
the case of $\dim \Pspan = 2$  and we denote by $C'$ the dimension 
$1$ component of $\Pspan \cap S$. The $\PP^3$'s containing 
$\Pspan$ generate a pencil of 
space curves $|D'|$ residual to $C'$. Notice that since plane curves 
in $S$ have degree at most $5$ we have $\deg D' \ge 6$.

\begin{prop} \xlabel{pDprime}
Let $D' \subset \PP^3$ be an irreducible space curve of degree at 
least $6$ that is contained in no quadric but in  a $a_3$-dimensional 
space of cubics, with $a_3\geq 3$.
Then $D'$ is either a septic of arithmetic genus $5$ and $a_3=3$ or a 
sextic of arithmetic genus
$a_3-1$ with $a_3\le4$.
\end{prop}

\begin{proof}
Let $Z$ be the curve component of the intersection of cubics that 
contain $D'$.
Then $\deg Z\leq 7$:  Since $Z$ lies in the complete intersection of 
two cubics, $\deg Z'\le
9$. But degree $9$ is impossible since such a curve only lies on two 
cubics.
Degree $8$ is also impossible since $Z$ would be linked $(3,3)$ to a 
line. Any
curve linked $(3,3)$ to a line is contained in precisely a pencil of 
cubics.  This latter result has a geometric version:  If the line is 
reduced in the complete intersection, then it intersects the linked 
curve in a scheme of length $4$.   Thus by Bezouts theorem, the line 
is contained in every cubic that contains the linked curve.

Since $D' \subset Z$ is irreducible and the residual part $D' - Z$ is 
at most a line, $Z$ is reduced.
Now, by assumption there are at least a net of pencils of cubic 
surfaces that contain $Z$.  So we consider the curve $E$ linked to 
$Z$ in a general such pencil.  By Bertinis theorem $E$ is singular 
only in the singular part of $Z$. Therefore $E$ must be reduced.

If $\deg Z = 7$ then $E$ must be a plane conic or two skew lines.  If 
$E$ is two skew lines, then at least one of them, say $L$, is not 
contained in $Z$. The union of $Z$ and the other line, say $L'$, is 
linked $(3,3)$ to $L$. By the geometric property of linkage above, 
$L$ intersects $Z$ in a scheme of length 
$4$, so it must be contained in $Z$, contrary to the above. Therefore 
$E$ must be a plane conic and $Z$ has arithmetic genus $5$ by the 
liaison formula.  Furthermore, in this case $a_3=3$.
In particular, if $a_3>3$, then $\deg D' \le \deg Z<7$. This proves 
our claim if $\deg D'=7$.

If $\deg D' = 6$, $a_3\geq 3$ then either $D'=Z$ or $D' = Z+L$ where 
$L$ is an additional line. In the latter case  $Z = D'\cup L$ is 
reduced, and linked $(3,3)$ to a plane conic $E$ by the previous 
argument. In particular $D'\cup L$ has arithmetic genus $5$. If the 
line $L$ lies in the plane of $E$, then $D'$ is linked $(3,3)$ to a 
plane cubic, so it is contained in a quadric, contrary to the 
assumption.  If $L$ meets $M$ in a point, then $L\cup E$ has 
arithmetic genus $0$, and by liaison, $D'$ has arithmetic genus $3$ 
and lies in $4$ cubics.  If $L$ does not meet $E$, then  $D'\cup E$  
is linked to the line $L$, so as above, $L$ intersects $D'$ in a 
scheme of length $4$. Thus $D'$ must have arithmetic genus $2$ and 
$a_3=3$.

If, on the other hand, $D'=Z$, then $E$ is a reduced curve of degree 
$3$.  If $E$ is not connected, then it has a line component that must 
intersect $D'$ in a scheme of length $4$ as above.  So this line 
would be contained in $Z$, contrary to our assumption.  Therefore $E$ 
is connected.  It is a plane cubic curve or a space curve of 
arithmetic genus $0$.  In the first case $D'$ would be linked $(3,3)$ 
to a plane cubic, and hence lie in a quadric hypersurface, against 
our assumption.  In the second case, $D'$ has arithmetic genus $3$ 
and $a_3=4$.
\end{proof}

\begin{prop} \xlabel{pa1} In the above notation,
$a \le 1$
\end{prop}

\begin{proof}
If $a = 2$ we  have $\dim \Pspan \ge 2$ by Corollary \ref{cPspanMin}.

If $\Pspan$ is a plane, then we consider the pencil of residual space 
curves
$|D'|$ introduced above.  As in the proof of Lemma \ref{Dsquare} we 
may assume that the general member $D'$ is a reduced and irreducible 
curve.

By Corollary \ref{cPspanMin} (\ref{iDimT}) the projection $\pispan$ 
is surjective.
By Lemma \ref{lContainingT} all $a+2=4$ quartics in the ideal of $C' 
\cup D'$ also
contain $T = \Pspan $. By Proposition \ref{pCohomologyEleven} the 
curve $C' \cup D'$ 
lies on no cubic. This implies that $D'$ lies on $4$ cubics and no 
quadric. Therefore $D'$ has degree $6$ and arithmetic genus $3$
by Proposition \ref{pDprime}. In this case we must have $C'=C$, 
$\Pspan = P$ and $D'=D$ since
$S$ contains only one plane quintic by Proposition \ref{plane5}.
But from Lemma \ref{Dsquare}
we know that $D$ is in fact smooth and irreducible of genus at most 
$2$, so we get a contradiction.

If $\Pspan = \PP^3$, i.e. $\Pperp = \PP^1$ the projection $\pispan$ is 
surjective: 
Every fiber is the intersection of a plane with
$\Prow \times \Pcol$ which contains at least $3$ points, and these 
points could not lie in $\Pperp $ since the Segre 
variety has no $3$-secant lines. If $\Pspan = \PP^4$ and $\Pperp$ is
a point,
then $T$ is a threefold since $\Prow \times \Pcol \subset \PP^5$ is
not a cone. If $\Pperp$ is outside of 
$\Pcol \times \Prow$, we have $\deg T =3 $ otherwise $\deg T = 2$.
In total $T$ is a threefold of degree at most $3$ contradicting Lemma
\ref{lContainingT}
\end{proof}

From now on we may assume $a=0$ or $1$.
\begin{prop} \xlabel{pColC}
If $a=0$, then $b=0$.
\end{prop}

\begin{proof}
If $a=0$, then $B = (q,l_1,l_2)^t$ has to be one of the types in 
Proposition \ref{pSecants}.
None of them has linear exterior syzygies.
\end{proof}

\begin{prop}\xlabel{pGindep}
If $a=1$ and $c=(l_1,l_2)^t$ is a column of $\Gamma$, then the linear 
forms $l_1$ and $l_2$ are
independent.
\end{prop} 

\begin{proof}
If $c=(0,0)^t$, then the restriction of $c$ to any point $p$ vanishes, 
and 
$h^1(\sI_S(5)|_p) \ge 1$. This is impossible. If after coordinate 
changes $c=(e_0,0)^t$, then
the restriction to $e_0$ vanishes and we obtain a contradiction as 
before.
\end{proof}
 
\begin{prop} \xlabel{pab1}
If $a=1$ and $b \ge 1$ then $\Pspan = \PP^1$.
\end{prop}

\begin{proof}
Since $B\Gamma=0$ the rows of $B$ have to be syzygies of the 
transpose of any column of $\Gamma$.
The syzygies of $c^t=(l_1,l_2)$ with $l_1$ and $l_2$ independent are 
generated by
\[
\begin{pmatrix}
l_1 & 0 & l_2 \\
0 & l_2 & l_1 \\
\end{pmatrix}.
\]
So the linear forms in $B$ must all lie in the span of the $l_i$, 
i.e. $\Pspan \subset  \langle l_1,l_e \rangle$.
Since by Corollary \ref{cPspanMin} the dimension of $\Pspan$ is at 
least one in this case, the
proposition follows.
\end{proof}

\begin{cor} \xlabel{cable1}
If $a=1$ then $b\le 1$.
\end{cor}

\begin{proof}
If $\Pspan \not= \PP^1$, then $b=0$. If $\Pspan = \PP^1$ we consider 
the linear part $B_1$ of
$B$. Possible columns of $\Gamma$ must be among the syzygies of 
$B_1$. After coordinate changes
there are only two possiblities for $B_1$:
\begin{enumerate}
\item If $B_1 = \left(\begin{smallmatrix} e_0 & 0 \\  0 & e_1 
\end{smallmatrix}\right)$ then the syzygies of $B_1$
are generated by $\left(\begin{smallmatrix} e_0 & 0 \\ 0 & e_1 
\end{smallmatrix}\right)$, so $b\le 2$. If $b=2$ then
$\Gamma =  \left(\begin{smallmatrix} e_0 & 0 \\ 0 & e_1 
\end{smallmatrix}\right)$ and some columns do not have full rank. 
This is impossible by Proposition \ref{pGindep}.
\item If $B_1 = \left(\begin{smallmatrix} e_1 & 0 \\ e_0 & e_1 
\end{smallmatrix}\right)$ then the syzygies of $B_1$
are generated by  $\left(\begin{smallmatrix} e_1 & 0 \\ e_0 & e_1 
\end{smallmatrix}\right)$. As before we obtain $b\le1$
\end{enumerate}
\end{proof}

\begin{cor}
The regularity of $\sI_S$ is at most $7$.
\end{cor}

\begin{proof}
If $b=0$ then then the regularity of $\sI_S$ is $6$. If $b=1$ then 
the regularity of 
$\sI_S$ is $7$ since $C = (l_1,l_2)^t$ does not have any linear 
syzygies, when the
$l_i$ are linearily independent. 
\end{proof}

Recall the well known fact

\begin{lem} \xlabel{l3Quadrics}
Let $E \subset \PP^3$ be a non degenerate pure $1$-dimensional scheme 
lying on a $3$-dimensional
set of quadrics. Then $E$ is defined by the $2 \times 2$ minors of a 
$2 \times 3$ matrix of linear form, i.e. a curve of degree $3$ and 
arithmetic genus $0$.
\end{lem}

\begin{proof}
Let $X$ be the scheme cut out by the $3$ quadrics. First we prove 
that $X$ can not contain a surface. 
Assume it does, then this surface must be a plane $P'$ and the 
residual scheme in $X$ is cut out by $3$ independent linear forms. So 
$X$ is $P'$ union a point. $E$ then must be contained in $P$ 
contradicting our assumptions.

So a general quadric in the ideal of $E$ is irreducible. The only 
curve 
on an irreducible quadric cut out by
exactly $3$ quadrics is a curve linked to a line in two quadrics, 
i.e. a cubic curve as described in the statement of the lemma.
\end{proof}

\begin{prop} \xlabel{pPspan2}
$\dim \Pspan \not= 2$.
\end{prop}

\begin{proof}
Assume that $\Pspan$ is a plane. Note that in this case only $a=1$ is 
possible. Let $C'$ be the
curve component of $\Pspan \cap S$ and
consider the pencil of residual space curves $|D'|$ introduced 
above.  Again, as in the proof of Lemma \ref{Dsquare}, we  may assume 
that a general member $D'$ is reduced and irreducible.

Since we consider $H := \PP^3$ containing $\Pspan$, the matrix $B_1$
vanishes on $H$ and the ideal of $S \cap H = C' \cup D'$ contains  a
$3$-dimensional space of quartics by Proposition \ref{pCohomology}
and no cubics by Proposition \ref{pCohomologyEleven}.

The projection $\pispan$ is surjective, since the preimage of a point 
is a line in
$\PP^3$ intersecting $\Prow \times \Pcolone$ in at least $2$ points
that can not both be contained in the projection center $\Pperp = 
\PP^0$.
So $\Pspan$ is
contained in all $3$ quartics of $I_{C' \cup D'}$ and therefore
$D'$ lies on three cubics and no quadrics.
On the other hand, $S$ has no plane curves of degree more than $5$, 
so $\deg D' \geq 6$.

Since $D'$ is reduced and irreducible, Proposition \ref{pDprime} 
applies, and so either $D'$ has degree $7$ and arithmetic genus $5$, 
or $D'$ has degree $6$ and arithmetic genus $2$.

In the second case of $\deg C' =11 - 6=5$ and we have
$C' = C$, $\Pspan = P$ and $D'=D$. By Lemma \ref{Dsquare}, the
general curve $D$ is smooth and irreducible of genus $D^2=2$.
In particular the intersection of $S$ with $P$ is the union
of $C$ and a scheme of length $2$. Therefore the ideal of $S$ is
not generated by sextics. Hence $b \ge 1$ and $\Pspan = \PP^1$ by 
Proposition \ref{pab1}.

This leaves us with the first case, where $D'$ has degree $7$ and arithmetic genus $5$ and $\deg C' = 11-7=4$. 
In this case $(D')^2=3$, so the pencil of curves $|D'|$ has a baselocus of length three in the plane $\Pspan$.  
In addition $\Pspan$ intersects $S$ in the quartic curve $C'$, so if the baselocus of $|D'|$ is disjoint from $C'$, 
then we may conclude from Example \ref{eQuartic3}
that $h^1(I_{\Pspan\cap S}(4))=3.$  Now $(C')^2=0$, so $h^1(\sO_{C'}(C'+2H))=0$, therefore, by Proposition \ref{pPlaneH2} (2),
we may compute the cohomology $h^1(I_{S}(4)|_{\Pspan})=h^1(I_{\Pspan\cap S}(4))=3.$  But by Proposition \ref{pCohomologyEleven},
$h^1(I_{S}(4)|_{\Pspan})\leq a+1=2$, which is in contradiction with the above.
What remains is to show that  the baselocus of $|D'|$ is disjoint from $C'$ in $\Pspan$.

The proof is rather roundabout, and we  start by considering again the projection $p$. 
There are two cases, either the projection center is contained in $\Prow \times \Pcolone$ or not. 
If not, then the pencil of lines given by the columns of $B_1$ does not have a basepoint. 
Since each of these line is at least a $5$ secant line to $S$, the curve $C'$ must be the plane quintic contrary to the above.
If the projection center is contained in $\Prow \times \Pcolone$, then two lines are blown down by the projection. 
The image of one of them is a basepoint of the pencil  given by the columns of $B_1$, say $e_0$.
 The other one corresponds to a column of $B_1$ whose entries only span a $\PP^0$ say $e_1$.
 Therefore after a coordinate change we can assume
\[
B = \begin{pmatrix}
q_1 & q_2 \\
e_0 & 0 \\
e_2 & e_1 \\
\end{pmatrix}
\]
Now consider the plane $P$ of the plane quintic. 
In the second column of $B$ we know by Proposition \ref{pSecants} (\ref{iSixInfGeometry}) that $e_1$ lies in $P$.
Furthermore, by this proposition, we know that $q_2 \not=0$. 
Therefore a general linear combination of columns must be of type (\ref{iGeneralGeometry}) in Proposition \ref{pSecants}. 
 Without loss of generality we can assume this to be the first column.   
  
By Proposition \ref{pSecants} (\ref{iGeneralGeometry}) we know that the line through $e_0$ and $e_2$ lies in a $\PP^3$ that also contains $P$.
 Therefore we have at least one point of $P$ on this line. This point can not be $e_1$ since $e_0$, $e_1$ and $e_2$ span $\PP^2$. 
 Therefore $P$ intersects $\Pspan$ in at least two distinct points i.e. in a line. If $P=\Pspan$ we are again in the case above.
  If $P \not = \Pspan$ their span $H$ is a $\PP^3$. The intersection curve $E = H\cap S$ then lies on $3$ quartics and no cubics. 
  The quartics contain $\Pspan$ as argued above and $P$, since $P$ contains a plane quintic. 
 Finally we have two possibilities; either the line of intersection $L=P\cap \Pspan$ is contained in $S$, or it is not.

 In the latter case, the residual curve $E-C-C'$ has degree $\deg(E-C-C') = 11-5-4=2$ and must lie on $3$ quadrics and no plane. This contradicts Lemma \ref{l3Quadrics}.

 \begin{lem} If $L=P\cap \Pspan$ is contained in $S$, then the baselocus of $|D'|$ in $\Pspan$ is disjoint from $C'$.
 \end{lem}
 \begin{proof}  The proof follows from a careful analysis of the reducible hyperplane section $E$.
 
 Notice first that the second column of $B$ vanishes on $L$, so,  
 by Proposition \ref{pSecants} (\ref{iSixInfGeometry}),  the selfintersection $L^2\leq -4$.
    The doubling of $L$ on $S$ is not contained in any of the planes $P$ or $\Pspan$, since such a doubling would mean that $L^2=1$ on $S$.
  Therefore $A=C'-L$ is a plane cubic with no component along $L$, and $B=C-L$ is a plane quartic curve with no component along $L$.  
  The two curves $A$ and $B$ can only intersect on $L$, but if they do, then this intersection point would be a singular point on $S$, a contradiction.
  Therefore $A\cdot B=0$.  On the other hand, $A\cdot L=3$ and $B\cdot L=4$, so the arithmetic genus  $p(A+B+L)=p(A+L)+p(B)+(A+L)\cdot B-1=9$.  
  The residual curve $G=E-A-B-L$ lies on $3$ quadrics and in no plane, so by Lemma \ref{l3Quadrics}
is has degree $3$ and arithmetic genus $p(G)=0$.
  
   The union $G+B$ lies in the pencil $|D'|= |H-C'|$.  
  Therefore $p(D')=p(G+B)=p(G)+p(B)+G\cdot B -1=5$, while $p(E)=p(G)+p(A+B+L)+G\cdot (A+B+L)-1=11.$
  Combined we get the intersection numbers $G\cdot B=3$ and $G\cdot (A+L)=0$.    
  
  If $G$ has a component along $L$, then $A+B+2L$ is contained in the union of the two planes $P$ and $\Pspan$, and $G-L$
   has degree two and is contained in at least 4 quadrics.  This means that $E$ is contained in $4$ quartics, contrary to the above.
   
   If $G$ has a component in common with $B$, then this component must have degree one or two. 
    In the latter case, the residual part of $G$ would be a line that intersects the first component.  
    But then it could not be a component of $A$, since $A$ and $B$ cannot intersect.  So
    $0=(A+L)\cdot G\geq L\cdot G= 2$, which is absurd. 
    
    The former case is similar if $G$ does not have a component in $A$.  
     If $G$ has a component in $A$, then it must be a line $L_{A}$, that does not intersect the line component $L_{B}$ of $G$ in $B$.
     Thus $G=L_{A}+L_{0}+L_{B}$, and $3=G\cdot B=(L_{0}+L_{B})\cdot B\leq 1+L_{B}\cdot B$, so $L_{B}\cdot B\geq 2 $.  
     But $(L_{B})^2= L_{B}\cdot B-L_{B}\cdot (B-L_{B})\geq 2-3= -1$.  On the other hand the intersection
      $(H-2L_{B})\cap L_{B}=(A+L+(B-L_{B})+L_{A}+L_{0})\cap L_{B}=(L+(B-L_{B})+L_{0})\cap L_{B}$ is finite of length at least five,
  while the intersection number $(H-2L_{B})\cdot L_{B}\leq 3$, a contradiction.

    Therefore $G$ has no component in $P$ and intersect $B$ and $P$ in a scheme of length $3$.  In particular $G$ does not intersect $L$,
     so $G$ has no component in common with $A$ and, since $A\cdot G=0$,  does not intersect $A$.
    In conclusion $G$ intersect $\Pspan$ in a scheme of length three outside $A\cup L$.
    \end{proof}


The lemma concludes the proof of the proposition.
 
\end{proof}

\begin{prop} \xlabel{pAssumeP1}
If $a=1$ and $\Pspan=L$ is a line, then $b=1$.
\end{prop}

\begin{proof}
We treat the possible matrices $B_1$ case by case:

If  $B_1 = \left(\begin{smallmatrix} e_0 & 0 \\ 0 & e_1 
\end{smallmatrix}\right)$, then
by Proposition \ref{pSecants} the plane $P$ contains $e_0$ and $e_1$. 
Furthermore any
line in $P$ though one of these points is either contained in $S$ or 
intersects $S$ in a scheme of length at least $6$. Since $P \cap S$ 
contains a plane quintic 
$C$ but no plane sextic $e_0$ and $e_1$ must be residual to $C$.

If  $B_1 = \left(\begin{smallmatrix} e_1 & 0 \\ e_0 & e_1 
\end{smallmatrix}\right)$, then 
by Proposition \ref{pSecants} the point $e_1$ is residual to $C$ in 
$P$. There are two cases
\begin{enumerate}
\item $e_0 \in P$. Then $B_1$ vanishes on every $\PP^3$ containing 
$P$. As in the proof of
Proposition \ref{pPspan2} the general residual space curve $D$ must
be of  degree $6$ and genus $2$. This implies that the residual
scheme in $P$ has length $2$.
\item $e_0 \not\in P$. In this case we can assume $P = e_1 \wedge e_2
\wedge e_3$ and
\[
B = \begin{pmatrix}
e_3 \wedge e_2 & e_1 & e_0 \\
e_3 \wedge e_2 + e_1 \wedge *& 0 & e_1
\end{pmatrix}^T
\]
by Proposition $\ref{pSecants}$. Furthermore $B$ must have a syzygy
$a = (e_1,q_1,q_2)$
with $q_1$ and $q_2$ forms of degree $2$. Equivalently
\[
B' = \begin{pmatrix}
e_3 \wedge e_2\wedge e_1 & e_1 & e_0 \\
e_3 \wedge e_2\wedge e_1& 0 & e_1
\end{pmatrix}^T
\]
must have a syzygy $(1,q_1,q_2)$ which is impossible since the first
row of $B'$ is independent
of the last two.
\end{enumerate}
In both case we therefore have a length $2$ scheme $R$ residual to 
$C$ in $P$.  Therefore $S$ is not cut out by sextics and $b=1$.

\end{proof}

\begin{rem} \xlabel{rP1}
If $a=1$ and $b=1$, then $L=\Pspan$ is a line in $P$ that contains a 
length two subscheme residual to $C$. The line $L$ is either a 
$7$-secant or $L \subset S$ with $L^2 \le -5$.
\end{rem}

We may summarize our results
\begin{cor}
For smooth surfaces in $\PP^4$ with Hilbert polynomial 
$P_S(n)=11/2n^2-9/2n+1$ there are three possible Hilbert functions, 
distinguished by their $B$ matrices. 
They differ from the Hilbert polynomials only in the degrees 
$n=1,\ldots ,5$ where they take the values\begin{enumerate}
\item $(5,15,35,70,116)$ and $B = (e_3 \wedge e_2, e_1, e_0)^T$, $B = 
(e_2 \wedge e_1, e_0,0)^T$,
or $B = (0,e_1,e_0)^T$;
\item $(5,15,35,69,116)$ and  $B = \left(\begin{smallmatrix} q_1& e_0 
& e_3 \\ q_2 & e_1 & e_4 
\end{smallmatrix}\right)^T$;
\item $(5,15,35,69,115)$ and $B = \bigl(\begin{smallmatrix} q_1& e_0 
& 0 \\ q_2 & 0 & e_1
\end{smallmatrix}\bigr)^T$ or $B = \left(\begin{smallmatrix} q_1& e_1 
& e_0 \\ q_2 & 0 & e_1 
\end{smallmatrix}\right)^T$.
\end{enumerate}
\end{cor}

\begin{proof}
The Hilbert function $H_S(n)$ is, of course, computed by 
$h^0(\sO_{\PP^4}(n)-h^0(\sI_S(n))$ so the difference to the Hilbert 
polynomial is the dimensions of the higher cohomology groups of 
$\sI_S(n)$.  The different possibilities for these groups are 
determined in Propositions \ref{pCohomologyEleven}, \ref{pa1}, 
\ref{pColC} and Corollary \ref{cable1},  and correspond to the 
following values of $a=h^1(\sI_S(4))-1$ and $b=h^1(\sI_S(5))$, 
namely  $a=b=0$, $a=1$ and $b=0$, and $a=b=1$.  The values of 
$h^0(\sI_S(n))$ and consequently of the Hilbert function then follows 
from the diagram of Proposition \ref{pCohomologyEleven}.
\end{proof}

Surfaces with the first Hilbert function were treated by Popescu in 
\cite{PopescuThesis}.
Examples of surfaces with the second Hilbert function were found by 
v. Bothmer, Erdenberger and Ludwig \cite {newfamily}. We treat this 
case in more detail in the next sections, both giving two different 
new constructions of such surfaces, and showing that they belong to a 
unique irreducible and unirational family.  In the third case there are
irreducible surfaces with this Hilbert function that are singular 
along $\Pspan$.  We have not found any smooth surface with this 
Hilbert function, but also cannot rule it out.

\section{Rational Surfaces with $d=11$, $\pi=11$, $a=1$, $b=0$} \xlabel{sa1b0}
\nosubsections 

Recall that, in the notation of the previous section, when $a=1$ and 
$b=0$ we have that $B_1$ is a $2 \times 2$ matrix and its
entries span a $\PP^3$.

\begin{prop}
Let $\Pspan = \PP^3$ and consider the smooth quadric $Q := \Prow 
\times \Pcolone \subset 
\Pspan$. Then the 
generic line $\Prow \times \{c\}$ in the pencil parametrized by the 
columns of $B_1$ intersects $S$ in a scheme of length at least $5$ 
and the special plane $P$ in one point. Furthermore $D^2<2$.
\end{prop}

\begin{proof}
Since $Q$ can not be contained in $S$ the generic 
line must intersect $S$ in a scheme of finite length. By Proposition 
\ref{pSecants} this length
is at least $5$. Since $Q$ is smooth, $P$ can not be contained in $Q$ 
and the intersection of the generic line with $P$ must be proper.  
Finally, $D^2\leq 2$ by Lemma \ref{Dsquare}.  But equality here means 
that $S$ intersects $P$ in a scheme of length two residual to $C$, so 
S is not cut out by sextics and $b>0$.
\end{proof}

We denote $C_{SP}$, $C_{SQ}$ and $C_{PQ}$ the curve components of $S 
\cap P$, $S \cap Q$ and $P \cap Q$ respectively. Also we set $H = 
\Pspan \cap S$ and $E = H-C_{SQ}$.

\begin{prop} $C_{SQ}$ is of type $(5,3)$ on $Q$ and $P$ intersects 
$Q$ in a line of type $(1,0)$.
\end{prop}

\begin{proof}
Since $B_1$ vanishes on $\Pspan$ the curve $H = S \cap \Pspan$ lies on $3$ 
quartics and
no cubic by Proposition \ref{pCohomology}. Since the general line in 
one of the rulings of $Q$ is at least a $5$-secant all quartics must 
contain $Q$. Therefore the residual curve $E = H - C_{SQ}$ lies on 
exactly $3$ quadrics and no cubic. By Lemma \ref{l3Quadrics} we 
obtain that $E$ is determinantal of degree $3$ and arithmetic genus 
$0$. Consequently $C_{SQ}$ has degree $8=11-3$.
Now the number of quintics independent of the quartics is at least 
\[
h^0(\sO_{\PP^3}(5))-h^0(\sO_{C_{SQ}}(5)) - 3\cdot 4 + s = 56 - 
(55+1-11) - 12 + s = s-1
\]
where $s$ is the number of linear syzygies between the quartics. 
Since the ideal of $D$ has $2$ linear syzygies between the $3$ 
quadrics, so do 
the three quartics obtained by multiplying with $Q$.
Therefore the ideal of $H$ contains at least one independent quintic. 
This shows that $C_{SQ}$ is
contained in a divisor of type $(5,5)$. So $C_{SQ}$ is of type 
$(5,3)$ or $(4,4)$. Since every line of type $(0,1)$ on $Q$ is at 
least a $5$-secant, the second choice is not possible.

Now consider the intersection $P \cap \Pspan$. If $P$ is a subset of 
$\Pspan$, $H$ must contain a plane quintic. Since $Q$ is smooth, $Q 
\cap P$ must be a plane conic in this case and consequently $E$ a 
plane cubic. This contradicts Lemma \ref{l3Quadrics}. So $L = P \cap
\Pspan$ is a line. By Proposition \ref{pSecants} it must intersect 
all $5$-secants on $Q$. Since the general divisor of type $(0,1)$ 
must be such a $5$-secant by the arguments above, we obtain that $L$ 
is of type $(1,0)$.
\end{proof}

Recall that reducible surface is called {\sl Zappatic} if
its components and the pairwise intersections
of two components are smooth. \cite{CilibertoMirandaZappatic}.

\begin{prop}
    The union $S\cup Q\cup P$ is 
    linked $(4,5)$ to a surface $B$ of degree $d=6$
    sectional genus $\pi = 3$. If $S \cup Q \cup P$ is Zappatic, then 
    $B$ is locally Cohen Macaulay with $\chi_B = 1$.
    \end{prop}

    \begin{proof}  First we note that any quintic hypersurface that 
    contains $S\cup P$ also contains $Q$. The reason is simply that 
    the intersection $(S\cup P)\cap Q = C_{SQ} \cup C_{PQ}$ is a 
curve of type $(6,3)$ on 
    the quadric.  Since the intersection $C_{SP}$ is a quintic 
curve, 
    the space of quintics in the ideal of $S$ that contains 
    $S\cup P\cup Q$ has codimension at most one in the space of all 
    quintics in this ideal.   Since $S$ lies on a unique quartic and 
on $10$ 
    quintics, this means that there is at least a $4$-dimensional 
    space of quintics that contains $S\cup P\cup Q$ independant of 
the quartic.

    Thus in the intersection of the unique quartic and a general 
    quintic in its ideal, the surface $S\cup P\cup Q$ is linked to a 
    surface $B$ of degree $6$.  

   The arithmetic genus $\pi_{U}$ of the union $U=S\cup P\cup Q$ is 
computed by 
the formula
\[
\pi_{U}=\pi_{S}+\pi_{P}+\pi_{Q}+deg(C_{SP})+deg(C_{SQ})+deg(C_{PQ})-2=11+5+8+1-2=23
\]   
The arithmetic genus $\pi_{B}$ of $B$ is computed by the formula 
    for liaison \cite{PS} of space curves:
    $$\pi_{U}-\pi_{B}=\frac{4+5-4}2(deg (U)-deg(B))=\frac {5\cdot 
    (14-6)}2=20,$$ 
    so $\pi_{B}=3$.

Now assume that $U$ is Zappatic. In particular we assume that the 
three curves
$C_{SP}$, $C_{SQ}$ and $C_{PQ}$ are smooth. Furthermore we assume
that the intersection $S\cap P \cap Q$ consists of $5$ distinct 
points  on the 
line of 
    intersction between $P$ and $Q$.  In three of these points all 
    three components meet pairwise in a curve, so the tangent cone of 
the union is 
    three planes that meet pairwise in a line. These points are 
Zappatic singularities
    of type $E_3$ in the notation of \cite{CilibertoMirandaZappatic}.
    
    In the remaining two points, the 
    isolated intersection points of $S\cap Q$, the plane $P$ 
    intersect both $S$ and $Q$ along a curve.  Therefore the tangent 
    cone to the union of the three surfaces at each of these points 
    is the cone over three lines that form a cubic space curve of 
    arithmetic genus $0$. So these points are Zappatic 
singularities
    of type $R_3$. Let $f=3$ be the number of $E_3$ singularities of 
the Zappatic
    surface $U$. Then 
    $$
    \chi_U 
    = \chi_S + \chi_Q + \chi_P - \chi_{C_{SP}} - \chi(C_{SQ}) - 
\chi(C_{PQ}) + f
    = 1 + 1 +1 + 5 + 7 - 1 + 3 = 17
   $$
   by the formula (3.16) of \cite{CilibertoMirandaZappatic}. Now we 
are in a position
   to compute $\chi_B$. We already computed $d_B = 6$ and $\pi_B = 
3$, so
   $\chi(\sO(dH)) = 3d^2+d+\chi_B$.    
   By the liason exact
   sequence
   $$
0 \to \omega_U \to \omega_{U \cup B} \to \sO_B(4) \to 0
  $$
  we have $\chi_{U \cup B} - \chi_U = 3\cdot 4^2+4+\chi_B$. But $U 
\cup B$ is a 
  complete intersection of type $(4,5)$ and therefore has $\chi_{U 
\cup B} = 70$.  
  Thus $\chi_B = 1$.
  
  Finally a Zappatic surface is locally Cohen Macaulay so by linkage 
$B$ is also 
  locally Cohen Macaulay.
\end{proof}

\section{counting dimensions} \xlabel{sDimensions}
\nosubsections

\begin{prop}
Let $S \subset \PP^4$ be a smooth surface of degree $d$, sectional
genus $\pi$
and Euler characteristic $\chi_S$ and $\sN$ the normal bundle on $S$
in $\PP^4$. Then
\[
\chi(\sN) = d(10-d)+5(\pi-1)+2\chi_S
\]
\end{prop}

\begin{proof}
By Hirzebruch-Riemann-Roch we have
\[
\chi(\sN) = \frac{1}{2}(c_1(\sN)^2-2c_2(\sN)) -
\frac{1}{2}c_1(\sN)K_S + 2 \chi_S
\]
on a surface. Substituting $c_1(\sN) = 5H+K$ and $c_2(\sN) = d^2$ and
applying the adjunction
formula $2\pi-2 = H(H-K)$, we obtain the formula above.
\end{proof}

\begin{rem}
Since $h^2(\sN)=0$ for a rational surface, $\chi(\sN)$ is the expected
dimension of the component of the
Hilbert scheme on which $S$ lies. \cite{GrothendieckHilbert}. 
In case $S$ is rational of degree $11$ and sectional genus
$11$,
we obtain
\[
\chi(\sN) = -11 + 5\cdot 10 +2 \cdot 1 = 41.
\]
\end{rem}
In the previous sections we analyzed smooth rational surfaces of 
degree
$11$ and sectional genus $11$ that lie on a quartic hypersurface and 
whose
ideal is generated in degree $6$, i.e. has Hilbert function
$H_S(n)=\{1,5,15,35,69, 11/2n^2-9/2n+1\}$.
More precisely we determined the possible linear parts of the maps of 
the Tate resolution $T(\sI_S)$
\[
\cdots\to 3E(-1) \xrightarrow{\alpha} E(-2)\oplus 2E(-3) \oplus E(-4)
\xrightarrow{\beta} 2E(-4)\to\cdots
\]

of such surfaces.

The entries of these maps are forms on $(\PP^4)^*$, so
the pure forms correspond to linear subspaces of $\PP^4$.  In 
particular,
the entries of the linear parts of the maps are points in $\PP^4$.  We
determined that the linear part of $\alpha$ has three entries that 
span a
plane $P=\PP^2$, while the linear part of $\beta$ is a $2\times 
2$-matrix
whose entries span $\Pspan=\PP^3$ and where the columns (and rows) 
span
the vertical (and horizontal lines) of a smooth quadric surface $Q$.
In this notation we count the parameters for the corresponding Tate
resolutions and find that they coincide with the expected dimension.
\begin{prop}
The family $F$ of minimal complexes
\[
3E(-1) \xrightarrow{\alpha} E(-2)\oplus 2E(-3) \oplus E(-4)
\xrightarrow{\beta} 2E(-4),
\]
such that the linear part of $\beta$ determine a smooth quadric 
surface
$Q\subset\Pspan = \PP^3$ and the linear part of $\alpha$ spans $P = 
\PP^2$
with $P \cap Q$ a line defined by a row of $B_1$,
is irreducible, rational and
of dimension $41$.
\end{prop}

\begin{proof}
First we can choose a plane $P \subset \PP^4$ and a smooth quadratic
surface $Q\subset\PP^3 \subset \PP^4$
such that $L= P \cap \PP^3=P\cap Q$ is a line.
For each such choice we can represent $\alpha$ and $\beta$ by 
matrices of
the form
\[
A = \begin{pmatrix}
e_0 & q_1 & q_4 & r_1 \\
e_1 & q_2 & q_5 & r_2 \\
e_2 & q_3 & q_6 & r_3 \\
\end{pmatrix}
\quad \quad
B = \begin{pmatrix}
p_1 & p_2 \\
e_0 & e_1 \\
e_3 & e_4 \\
0 & 0 \\
\end{pmatrix}
\]
where $P=\langle e_0, e_1, e_2\rangle$, $L=\langle e_0, e_1\rangle$ 
and $\Pspan=\langle e_0,e_1,e_3,e_4\rangle$.
The entries  $p_i$ and $q_i$ are $2$-forms and $r_i$ are $3$-forms. 
Notice
that the
relation $AB = 0$ is linear in the coefficients of the $2$-forms.
In the corresponding linear
system, we have $80$ coefficients of $2$-forms in $A$ and $B$ and six
$3$-forms with $10$ coefficients in $AB$, so we expect a 
$20$-dimensional
solution.  But in fact the relations of $AB=0$ are dependent and we 
find a
$25$-dimensional affine solution space. See \cite{linesplanesweb} for 
the calculation.
Now the
quadratic part of both matrices is only defined modulo the linear 
part.
Projectively we obtain a $14 = 25 - 2\cdot5 -1$ dimensional
solution space. The cubics can be chosen arbitrarily, but are defined
only modulo the linear and quadratic part of $A$.  Furthermore the two
syzygies of $A$ given by $B$ are degree $3$ dependencies between the
linear and the quadratic part, so we obtain a $11 = 30 - 1\cdot 10 - 2
\cdot 5 -1 +2$ dimensional space of possible degree $3$ parts of $A$.
In total we have shown that the family of complexes $F$ is 
birationally
parameterized by an irreducible and Zariski open set
\[
X \subset \GG(2,5) \times (\PP^4)^*
\times \PP_{\text{quadrics}}^6
\times \PP_{\text{deg2part}}^{14}
\times \PP_{\text{deg3part}}^{11}.
\]
 In particular $\dim X = 6+4+6+14+11 =
41$.
\end{proof}

\begin{thm}\xlabel{tcons}
The family of smooth rational surfaces of degree $11$, sectional genus
$11$ with Hilbert function $H_S(n)=\{1,5,15,35,69, 11/2n^2-9/2n+1\}$ 
is
unirational, irreducible of dimension $41$. The general member of the 
family
is linked $(4,5)$ to a Zappatic surface $P\cup Q\cup B$, where $P$ is 
a
plane, $Q$ is a smooth quadric surface and $B$ is a smooth Bordiga
surface, such that $L=P\cap Q$ is a line, $B\cap Q$ consists of three 
distinct
lines that intersect $L$, and $P\cap B$ is the union of a line 
distinct
from $L$ and two points on $L$.  In particular $S\cup P\cup Q$ is an
arithmetically Cohen-Macaulay surface defined by the $4\times 4$ 
minors of
a $4\times 5$ matrix with $4$ columns of linear forms and $1$ column 
of
quadratic forms.
\end{thm}
\begin{proof}
By the previous proposition it remains for the first part to give an
example.  This is done by choosing random matrices $A$ and $B$ 
satisfying
the above conditions, computing the minimal free resolution of 
$\sI_S$ via
the BGG-Correspondence and checking that this is an ideal sheaf 
defining a smooth surface by the
Jacobian criterion.  An effective procedure is to compute an example 
over
a finite characteristic.  This is done with with Macaulay 2 and 
documented at \cite{linesplanesweb}.
The fact that the general member is linked $(4,5)$ to a Zappatic 
surface
is an open condition that is also checked in an example.
Finally, since a Bordiga surface is an arithmetically Cohen-Macaulay
surface defined by the $3\times 3$ minors of a $3\times 4$ matrix with
linear entries, the linked surface $S\cup P\cup Q$ is also 
arithmetically
Cohen-Macaulay, and the proposition follows.
\end{proof}

\begin{rem}
The example by  v.\,Bothmer, Erdenberger and Ludwig is a blowup 
of the plane in $20$ points.  In fact the linear system has 
the form
$9L-3E_{1}-\sum_{i=2}^{15}2E_{i}-\sum_{i=16}^{20}E_{i}$
where $L$ is the pullback of a line from the plane, while the $E_{i}, 
i=1,\ldots,20$ are the exceptional divisors of the blowup. 
From Theorem \ref{tcons} it follows that the linear system of any 
smooth surface 
with the second Hilbert function has this form.
The challenge remains to determine necessary and sufficient conditions
for position of the 20 points in $\PP^2$.
\end{rem}

\section{construction} \xlabel{sConstruction}
\nosubsections

We can use the properties of a surface $S$ as in the Theorem 
\ref{tcons} to
give a geometric construction

\begin{construction} $\quad$ \xlabel{cSurface}

\begin{enumerate}
\item Choose a line $L$ in $\PP^2$
\item Choose general points $P_1$, $P_2$ and $P_3$ on $L$
\item Choose general points $P_4,\dots,P_8$ outside of $L$
\item Let $C$ be the unique irreducible quartic curve that contains
$P_1,\dots,P_8$ and is singular in $P_3$, $P_4$ and $P_5$
\xlabel{iUnique4}
\item Choose general points $P_9$ and $P_{10}$ on $C$
\item Let $B$ be the blowup of $\PP^2$ in $P_1,\dots, P_{10}$ and
denote
the exceptional divisors by $E_1,\dots, E_{10}$.
\item Embed $B$ in $\PP^4$ with the linear system $|H| :=
|4L-E_1-\dots-E_{10}|$ the image will be
a Bordiga surface of degree $6$ and sectional genus $3$
\item Let $\Pspan$ be the hyperplane in $\PP^4$ corresponding to $C
\in |H|$. Since $C$ is singular in $P_3$, $P_4$ and $P_5$, 
$\Pspan$ contains the exceptional lines $E_3$, $E_4$ and $E_5$.
\item Let $Q \subset \PP^3_C$ be the unique quadric containing these
lines.
\item Let $\widetilde{L}$ be the strict transform of $L$. It is again
a line in $\PP^4$. Let
$P$ be the unique $\PP^2$ containing $\widetilde{L}$ and intersecting
$E_4$ and $E_5$.
\item Let $S$ be a $(4,5)$ linkage of $P \cup Q \cup B$.
\end{enumerate}
\end{construction}

\begin{thm}\xlabel{tlinkageconstruction}
The construction \ref{cSurface} yields a $41$-dimensional unirational, irreducible
family of smooth rational surfaces $S$ of degree $11$ and sectional
genus $11$ with precisely two $6$-secants lying on a unique quartic.
\end{thm}

\begin{proof}

It is straightforward to check that each step of the construction is
possible except the last one.
For the last step we need to show that $P\cup Q\cup B$ lies on a 
quartic
and a quintic hypersurface with no common component. For this we first
consider the exact sequence of ideal sheaves
\[
0\rightarrow I_B(2)\rightarrow I_{B\cup Q}(3) \rightarrow I_{B\cup 
Q}|_H(3)\rightarrow 0
\]
where the first map is multiplication by the linear form defining the
hyperplane $H$ that contains $Q$.  All cohomology groups on the left
and on the right vanish:  On the left $B$ is arithmetically Cohen 
Macaulay and does not lie on any quadric, while on the right 
$I_{B\cup Q}|_H$ is the ideal sheaf in $\sO_H$ of the union of $Q$ 
and the twisted cubic curve $C$ on $B\cap H$ residual to $Q$.  
Therefore the cohomology of $I_{B\cup Q}|_H(3)$ coincides with the 
cohomology of $I_C|_H(1)$, which vanishes.  We conclude that the 
cohomology groups of the sheaf in the
middle also vanishes.  Similarly, twisting with $\sO_{\PP^4}(1)$, we 
get
$h^0(I_{B\cup Q}(4))=7$ and $h^1(I_{B\cup Q}(4))=0$.
Next, consider the exact sequence of ideal sheaves
\[
0\rightarrow I_{B\cup Q}(3)\rightarrow I_{B\cup Q\cup P}(4) 
\rightarrow
I_{B\cup Q\cup P}|_{H'}(4)\rightarrow 0
\]
where the first map is the multiplication by the linear form defining 
a
general hyperplane $H'$ through $P$.  By the above, $B\cup Q\cup P$ 
lies in
a quartic hypersurface if and only if $(B\cup Q\cup P)\cap H'$ lies 
on a
quartic surface in $H'$.  But $(B\cup Q\cup P)\cap H'$ is the union 
$P\cup
E\cup L'$, of the plane $P$, an elliptic quintic curve $E$ in $H'\cap 
B$
with a trisecant line $L$ in $P$ and a line $L'$ such that $L\cup 
L'=Q\cap
H'$. So $(B\cup Q\cup P)\cap H'$ is contained in a quartic surface if 
and
only if $E\cup L'$ is contained in a cubic surface.  Since $E\cup 
L=H'\cap
B$ lies in $4$ cubics, and $L'$ meets $L$, there is at least one cubic
surface that contains $E\cup L'$, and hence at least one quartic
hypersurface that contains $B\cup Q\cup P$.

On the other hand, there are $15$ quartics that contain $E$, i.e. $10$
quartics that contain $E\cup L'$. Hence, as above,
\[
h^0(I_{B\cup Q\cup P}(5))=h^0(I_{B\cup Q}(4))+ h^0(I_{B\cup Q\cup 
P}|_{H'}(5)) =7+10=17
\]
and $B\cup Q\cup P$ lies in $12$ quintic hypersurfaces that are
independent of the quartic.  Since the Bordiga surface is not 
contained in
any reducible cubic hypersurface, $B\cup Q$ is not contained in any 
cubic
hypersurface and $Q\cup P$ is not contained in any hyperplane, any 
quartic
that contains $B\cup Q\cup P$ must be irreducible.  Therefore the 
general
quintic and quartic that contains $B\cup Q\cup B$ have no common 
component
and so
$B\cup Q\cup B$ is linked $(4,5)$ to a surface $S$.

The parameters involved in the construction form
an open set in
\[
Y\subset(\PP^2)^* \times (\PP^1)^3 \times (\PP^2)^5 \times (C)^2 
\times
\Aut(\PP^4) \times
\PP^{11}_{linkage}.
\]
 Since $C$ is rational
this proves that our family is unirational.
To find the dimension of the Hilbert scheme component we need to 
subtract dimension of the automorphisms of $\PP^2$ (projective 
dimension 8) and
the dimension of the space of independent quintic hypersurfaces that 
contain $S \cup
P \cup Q$ (projective dimension 3 by Theorem \ref{tcons}) , i.e. the 
space of Bordiga surfaces that lead to the same
$S$.  Therefore these surfaces $S$ belong to a

\[
2+3+10+2+24+11-8-3 = 41
\]
dimensional family in the Hilbert scheme.

Next we compute the Hilbert polynomial of $S$.

Of course, the degree of $S$ is $11$. The sectional genus is given by 
the
liaison formula
\[
\pi_S-\pi_{B\cup P\cup Q}=5/2(11-9)=5
\]
i.e. $\pi_S=11$ since the sectional genus of $B\cup Q\cup P$ is $6$.
To get the Euler characteristic $\chi_S$ we first compute $\chi_{B\cup
Q\cup P}$.  By construction $B\cup Q\cup P$ is Zappatic, i.e.
the three components $B$, $P$, $Q$ and their pairwise intersections 
are
smooth.  Furthermore the intersection $B\cap Q\cap P$ consists of
precisely three points on the line $L=P\cap Q$.  At one of them, the 
point
$L\cap L'$, the three surface compontent intersect pairwise in 
codimension
$1$, so this point is a Zappatic singularity of type $E_3$ on the 
union. 
At the two other points, $B$ and $P$ intersect in codimension $2$, 
while
the $B$ and $Q$ and $Q$ and $P$ intersect in codimension $1$, so these
points are Zappatic singularities of type $R_3$.
The Euler characteristic of $B\cup Q\cup P$ is therefore
\[
\chi_{B\cup Q\cup P}=\chi_B+\chi_Q+\chi_P-\chi_{L'}-\chi_L-\chi_{B\cap
Q}+f=1+1+1-1-1-3+f=f-2
\]
by the formula for Zappatic surfaces.  Since $f=1$ count the number of
Zappatic singularities of type $E_3$, we get $\chi_{B\cup Q\cup 
P}=-1$.
In the liaison exact sequence
\[
0\rightarrow \omega_{B\cup Q\cup P}\rightarrow 
\omega_{(4,5)}\rightarrow
\sO_S(4)\rightarrow 0
\]
The Euler characteristic of the first two sheaves are $-1$ and $70$
respectively, while the Hilbert polynomial of $\sO_S$ is $P_S(d)=11/2
d^2-9/2 d+\chi_S$.  In particular $P_S(4)=70+\chi_S$, so by the exact
sequence $\chi_S=1$.

To see that the general surface of this family is of the kind found 
in the
previous section, we have checked one example for smoothness (see 
\cite{linesplanesweb}). 
Since Popescu \cite{PopescuThesis} showed that there are no 
nonrational surfaces of these
invariants, $S$ must be rational.  Furthermore, by liaison, $S\cap Q$ 
is a
curve of type $(3,5)$, while $L$ is of type $(0,1)$.  So two of the 
three
lines in $B\cap Q$ are $6$-secants to $S$.  
\end{proof}


\def\cprime{$'$} \def\cprime{$'$}

\end{document}